\documentclass[11pt]{article}
\usepackage{amsfonts,amsmath,amsthm,amssymb,mathtools}
\usepackage{mathrsfs,array,xcolor,authblk,indentfirst}
\usepackage[numbers,sort&compress]{natbib}
\usepackage[a4paper,left=2.5cm,right=2.5cm,top=2cm,bottom=2cm]{geometry}
\usepackage[colorlinks,linkcolor=blue,citecolor=blue]{hyperref}

\allowdisplaybreaks

\newtheorem{theorem}{\indent Theorem}[section]
\newtheorem{lemma}[theorem]{\indent Lemma}
\newtheorem{proposition}[theorem]{\indent Proposition}
\newtheorem{corollary}[theorem]{\indent Corollary}

\newtheorem{remark}[theorem]{\indent Remark}
\numberwithin{equation}{section}
\newcommand{\R}{\mathbb R}
\newcommand{\N}{\mathbb N}
\newcommand{\Z}{\mathbb Z}
\newcommand{\T}{\mathbb T}
\newcommand{\E}{\mathbb E}
\newcommand{\Prob}{\mathbb P}

\newcommand{\cAop}{\mathcal A}

\newcommand{\cZ}{\mathcal Z}

\newcommand{\cB}{\mathcal B}
\newcommand{\Pt}{P_t}
\newcommand{\Osc}{\operatorname{Osc}}
\newcommand{\cP}{\mathcal P}
\newcommand{\cM}{\mathcal M}
\newcommand{\cV}{\mathcal V}

\newcommand{\cJ}{\mathcal J}
\newcommand{\cK}{\mathcal K}
\newcommand{\fB}{\mathfrak B}
\newcommand{\fc}{\mathfrak c}
\newcommand{\fm}{\mathfrak m}
\newcommand{\fw}{\mathfrak w}
\newcommand{\Lt}{L_t}
\newcommand{\Ezero}{E_0}
\newcommand{\dd}{\,\mathrm d}
\newcommand{\norm}[1]{\left\| #1 \right\|}

\newcommand{\ip}[2]{\left\langle #1,\, #2\right\rangle}
\newcommand{\Lam}{\Lambda}
\newcommand{\PV}[1]{P_{#1}^{V}}
\DeclareMathOperator{\supp}{supp}
\DeclareMathOperator{\spn}{span}

\begin{document}
\title{\textbf{Large deviation principle for occupation measures of the 2D stochastic Boussinesq equations with highly degenerate noise}}
\author[a]{Xue Wang\thanks{wangxuexue1997@163.com}}
\author[a]{Jianhua Huang\thanks{jhhuang32@nudt.edu.cn}}
\author[b]{Jiangwei Zhang\thanks{Corresponding author. zjwmath@163.com}}
\affil[a]{\small\itshape College of Science, National University of Defense Technology, 
		Changsha, Hunan 410073, P.R.China}
\affil[b]{\small\itshape Institute of Applied Physics and Computational Mathematics, 
    	Beijing, 100088, P.R.China}
\date{}
\maketitle

\begin{abstract}
	We study large deviations for occupation measures of the two-dimensional stochastic Boussinesq equations on the torus, driven by additive white noise acting only on four Fourier modes of the temperature equation.	
	First, we establish a multiplicative ergodic theorem for Feynman–Kac semigroups with bounded continuously differentiable potentials having bounded derivatives. The proof combines weighted energy estimates and approximate controllability with a uniform Feller estimate. For the latter, conditional Malliavin estimates give residual decay sufficient to offset the growth of the Feynman--Kac weights.
	Next, we prove a large deviation principle for occupation measures in the weak topology, uniformly over initial distributions with a common exponential moment bound, using pressure convergence, exponential tightness and uniqueness of equilibrium states for smooth potentials. In particular, the good rate function vanishes only at the unique invariant measure.
	Finally, a Wiener small-ball argument gives a probability measure with finite, strictly positive rate supported in each closed ball centered at the origin with radius in $(0,1]$.
	
\end{abstract}

{\hspace*{2mm}AMS Subject Classification (2020): 60F10, 60H15, 35Q35.}

\section{Introduction}\label{sec:intro}
In the study of randomly forced fluids, invariant measures and ergodic theorems describe the long-time statistics of the flow.
These results characterize typical time averages, but do not by themselves determine the probabilities of deviations from the limiting statistics. Large deviation principles quantify these probabilities, providing a framework for studying rare events in fluid and climate models~\cite{RWB18,GL21}.
For buoyancy-driven flows, the Boussinesq equations couple fluid motion and temperature, as in models of Rayleigh-B\'enard convection and atmospheric and oceanic circulation~\cite{Maj03}.
The aim of this article is to establish a large deviation principle for occupation measures when the noise acts on only four temperature modes.

We consider the following 2D stochastic Boussinesq equations on the torus $\T^2=\R^2/(2\pi\Z)^2$:
\begin{equation}\label{eq:bouss-intro}
\left\{
\begin{aligned}
&\dd u+(u\cdot\nabla)u\dd t=\big(-\nabla p+\nu_1\Delta u+g\theta e_2\big)\dd t,\\
&\dd\theta+(u\cdot\nabla)\theta\dd t=\nu_2\Delta\theta\dd t+\sigma_\theta\dd W,\\
&\nabla\cdot u=0,
\end{aligned}
\right.
\end{equation}
where $u=(u_1,u_2)$, $p$ and $\theta$ denote the velocity, pressure and temperature.
The constants $\nu_1,\nu_2>0$ are the viscosity and thermal diffusivity.
The buoyancy coefficient is $g\neq0$, and $e_2=(0,1)^\top$.
The noise acts on four temperature modes:
\[
\sigma_\theta\dd W=\alpha_1\cos x_1\dd W^1+\alpha_2\sin x_1\dd W^2+\alpha_3\cos x_2\dd W^3+\alpha_4\sin x_2\dd W^4,\qquad\alpha_i\neq 0.
\]
The velocity receives no direct stochastic forcing.
Randomness reaches it through buoyancy and is transferred to other modes by nonlinear transport.
Writing $w=\partial_1u_2-\partial_2u_1$, we use the state variable $U=(w,\theta)$ in $H=L^2_0\times L^2_0$.
The precise setting is given in Section~\ref{sec:setting}.
Our object of study is the occupation measure
\[
\Lt:=\frac1t\int_0^t\delta_{U_s}\dd s,\qquad t>0,
\]
which records time averages through $t^{-1}\int_0^tf(U_s)\dd s=\ip f{\Lt}$.

The Boussinesq equations have been studied from both analytical and probabilistic points of view.
For the 2D deterministic system, global well-posedness with only one of the two dissipation terms was established by Hou and Li~\cite{HL05} and Chae~\cite{Cha06}.
In the stochastic setting, small-noise large deviation principles were proved by Duan and Millet~\cite{DM09} for multiplicative Gaussian forcing and by Zheng and Huang~\cite{ZH16} for L\'evy forcing.
For the present degenerate forcing, F\"oldes et al.~\cite{FGRT15} proved unique ergodicity and mixing, together with a law of large numbers and a central limit theorem for time averages.
Related ergodic properties for randomly forced Rayleigh-B\'enard convection were obtained in~\cite{FGRW16}, as well as for the degenerately forced 2D Navier-Stokes, fractional magnetohydrodynamic and real Ginzburg-Landau equations~\cite{DP18,PHZ20a,PHZ20b}.
These ergodic results describe typical long-time statistics and motivate the study of atypical occupation measures as $t\to\infty$ at fixed noise intensity.

The study of large deviations for occupation measures goes back to the work of Donsker and Varadhan~\cite{DV75}.
For infinite-dimensional systems, Wu~\cite{Wu01} introduced a criterion based on the strong Feller property, irreducibility and hyper-exponential recurrence. Under suitable non-degeneracy assumptions on the noise, this criterion has been applied to the stochastic Navier-Stokes and Burgers equations~\cite{Gou07a,Gou07b} and to monotone SPDEs~\cite{WXX21}.
Related applications to dissipative fluid and reaction diffusion equations can be found in~\cite{KM23a,KM23b}. The approach developed in~\cite{JNPS15,JNPS18} for randomly kicked dissipative PDEs combines a Kifer-type criterion~\cite{Kif90} with multiplicative ergodicity of Feynman-Kac semigroups. A continuous time version of this approach was developed by Martirosyan and Nersesyan~\cite{MN18a}.
For Navier-Stokes equations driven by white noise, Nersesyan~\cite{Ner19} used this framework to establish a large deviation principle without assuming the strong Feller property.
Nersesyan et al.~\cite{NPX23} subsequently used Malliavin calculus to establish the large deviation principle of occupation measure under highly degenerate forcing satisfying the condition of Hairer and Mattingly~\cite{HM06}.
Their multiplicative ergodic theorem applies to general potentials in $C_b^1(H)$.
Liu and Lu~\cite{LL26} established a Donsker-Varadhan large deviation principle for 2D Navier-Stokes equations with a periodic force and degenerate noise, using time-inhomogeneous Feynman-Kac operators.
Related results for space-time localized forcing were obtained by Peng and Xu~\cite{PX26}.
A level-3 large deviation principle for the Navier-Stokes equations in bounded domains was established in~\cite{Zhao25}.

Controllability provides a further connection between the deterministic dynamics and the stochastic system under degenerate forcing.
For the Boussinesq equations, Glatt-Holtz et al.~\cite{GHM18} proved approximate controllability using the four temperature modes specified above, together with exact controllability on finite dimensional projections.
These properties hold at every positive time and imply full support of the stochastic transition probabilities.
Control by spatially localized temperature forcing was studied by Nersesyan and Rissel~\cite{NR25} and Rissel~\cite{Ris26}, while localized controls in both equations were considered in~\cite{Ris25}.
Related controllability results for the three-dimensional primitive equations were obtained in~\cite{BGN23}. 
For the system considered here, noise enters only through the temperature equation, and the velocity is driven through the buoyancy coupling.
We establish the occupation-measure large deviation principle for this coupled system at fixed noise intensity.

Our main result is a level-2 large deviation principle in the weak topology, uniformly over initial distributions with a common exponential moment bound (see Theorem~\ref{thm:A}).
We also prove multiplicative ergodicity of the Feynman-Kac semigroups for $V\in C_b^1(H)$, with convergence for exponentially growing terminal observables (see Theorem~\ref{thm:B}).
The good rate function is the Legendre transform of the pressure and has the invariant measure $\mu_*$ as its unique zero.
For every $\varepsilon\in(0,1]$, there is a probability measure supported in $B_H(0,\varepsilon)$ with finite positive rate.

We use the Feynman-Kac approach of~\cite{JNPS18,Ner19,NPX23}.
The Malliavin covariance estimates for the temperature-vorticity system are supplied by~\cite{FGRT15}, and approximate controllability by~\cite{GHM18}.
The estimates needed for exponential trajectory weights are developed in Sections~\ref{sec:apriori} and~\ref{sec:fk}.
First, a weighted energy controls the buoyancy coupling and yields exponential moments for the state and its integrated $H^1$ norm.
The resulting dissipative entrance bound, followed by approximate control from a fixed bounded set, gives an irreducibility probability independent of the initial radius.
The required control time may depend on that radius.

In the Feynman-Kac gradient estimate, the exponential trajectory weight can offset the decay of the control residual. To prove the uniform Feller property, we retain an exponential weight in the residual estimate at times $2n$.
Conditioning on the past then combines this estimate with bounds for the Malliavin derivatives of the restarted flow and for the Skorokhod integral of the control.
The regularization parameter is chosen after $V$, so that the residual decay compensates for the growth of the Feynman-Kac weight.
Following the integration-by-parts argument of~\cite{NPX23}, we obtain
\[
|\nabla_\xi P_t^V\psi(U)|\le C e^{\eta\norm U^2}\norm{P_t^V\mathbf1}_{R_0}
\big[(1+\norm{\nabla V}_\infty)\norm\psi_\infty+\varepsilon_*^t\norm{\nabla\psi}_\infty\big]\norm\xi .
\]
The exponential occupation bound extends the pressure limit to bounded continuous potentials, and a Doob transform gives uniqueness of equilibrium states for smooth potentials.
Finally, stability under uniformly small driving paths and the Wiener small-ball bound give finite-rate measures supported near the origin.

The rest of the paper is organized as follows.
Section~\ref{sec:setting} presents the setting and main results.
Section~\ref{sec:apriori} establishes the properties of the path solution map and the a priori estimates.
Section~\ref{sec:fk} is devoted to the properties of the Feynman-Kac semigroups, and the main results are proved in Section~\ref{sec:proofs}.
\section{Setting and main results}\label{sec:setting}

In this section, we first introduce the vorticity-temperature formulation and the path solution map. We then define the Feynman-Kac semigroups and their weights, and state the main results.

\subsection{Notation and abstract formula}\label{sec:notation}

We use $C$ to denote a positive constant whose value may change from one estimate to the next. Relevant parameter dependence is stated explicitly. Constants that must remain fixed within an argument are distinguished by subscripts.

Let $\N:=\{0,1,2,\dots\}$, $\T^2=\R^2/(2\pi\Z)^2$ and $\Z^2_0:=\Z^2\setminus\{0\}$.
 Denote by $L^2_0$ the real-valued mean-zero functions in $L^2(\T^2)$.
For scalar functions, $\ip\cdot\cdot$ and $\norm\cdot$ are the inner product and norm of $L^2(\T^2)$.
The vorticity is $w=\partial_1u_2-\partial_2u_1$.
The mean-zero velocity is recovered by the Biot-Savart law $u=K*w:=\nabla^\perp\Delta^{-1}w$, where $\nabla^\perp:=(-\partial_2,\partial_1)$.
In particular, $\nabla\cdot(K*w)=0$.
For a vorticity $c$ and a scalar function $f$, write $b(c,f):=(K*c)\cdot\nabla f$.
Let $\Z^2_+:=\{k\in\Z^2_0:k_2>0\}\cup\{(k_1,0):k_1>0\}$.
For $k\in\Z^2_+$, define the temperature directions $\sigma_k^0:=(0,\cos k\cdot x)$ and $\sigma_k^1:=(0,\sin k\cdot x)$.
The corresponding vorticity directions are $\psi_k^0:=(\cos k\cdot x,0)$ and $\psi_k^1:=(\sin k\cdot x,0)$.

The phase space is $H:=L^2_0\times L^2_0$. Let $\fc:=\nu_1\nu_2/g^2$. For $U=(w,\theta)$ and $\tilde U=(\tilde w,\tilde\theta)$ we define the weighted inner product, the weighted norm and the homogeneous weighted Sobolev norms
\begin{gather*}
\ip U{\tilde U}_{\fc}:=\fc\ip w{\tilde w}+\ip\theta{\tilde\theta},\qquad
\norm U^2:=\ip UU_{\fc},\\
\norm U_{H^s}^2:=\fc\norm{(-\Delta)^{s/2}w}^2+\norm{(-\Delta)^{s/2}\theta}^2\quad(s\ge0).
\end{gather*}
Let $H^s:=\{U\in H:\norm U_{H^s}<\infty\}$, with inner product $\ip\cdot\cdot_{H^s}$.
Thus $H^0=H$.
Unless specified otherwise, inner products, adjoints, orthogonal projections and operator norms on $H$ refer to $\ip\cdot\cdot_\fc$.
Scalar $L^2$ spaces and $\R^4$ retain their usual inner products.
On mean-zero functions, $\norm\cdot_{H^s}$ is equivalent to the inhomogeneous Sobolev norm, and $\norm U\le\norm U_{H^1}$.
We write $U_w$ and $U_\theta$ for the vorticity and temperature components of $U$.
Let
\begin{gather*}
B(U,\tilde U):=\big(b(w,\tilde w),\,b(w,\tilde\theta)\big),\qquad B(U):=B(U,U),\\
AU:=(-\nu_1\Delta w,-\nu_2\Delta\theta),\qquad GU:=(g\partial_1\theta,0),
\end{gather*}
and set $F(U):=-AU-B(U)+GU$, $\kappa:=\min\{\nu_1,\nu_2\}$ and $Z(U):=\ip{AU}{U}_{\fc}\ge\kappa\norm U_{H^1}^2$.

The noise acts only on the four lowest modes of the temperature equation. Let $\Ezero:=\spn\{\cos x_1,\sin x_1,\cos x_2,\sin x_2\}$. Throughout the paper, we assume the following:
\begin{itemize}
\item[(A1)] $\nu_1>0$, $\nu_2>0$, $g\neq0$.
\item[(A2)] $\cZ=\{(1,0),(0,1)\}$, $\alpha_k^l\neq0$ for all $k\in\cZ$ and $l\in\{0,1\}$, and
\[
\sigma_\theta\dd W=\sum_{k\in\cZ,\,l\in\{0,1\}}\alpha_k^l\sigma_k^l\dd W^{k,l},\qquad
\fB_0:=\norm{\sigma_\theta}_{HS}^2=\sum_{k\in\cZ,\,l\in\{0,1\}}(\alpha_k^l)^2\norm{\sigma_k^l}^2 .
\]
\item[(A3)] The stochastic basis $(\Omega,\mathcal F,\{\mathcal F_t\},\Prob)$ satisfies the usual conditions.
The process $W=(W^{k,l})_{k\in\cZ,\,l\in\{0,1\}}$ is a standard four-dimensional Brownian motion.
\end{itemize}
Above, $\sigma_\theta:\R^4\to\{0\}\times\Ezero\subset H$.
The coefficients and Brownian motions in the introduction are now indexed by $(k,l)\in\cZ\times\{0,1\}$.
Fix an ordering of these four indices and write $W^j:=W^{k,l}$ for $j=1,\dots,4$.
Let $\mathbf e_j$ denote the standard basis of $\R^4$.
The vorticity-temperature formulation of \eqref{eq:bouss-intro} is
\begin{equation}\label{eq:abstract}
\dd U+\big(AU+B(U)\big)\dd t=GU\dd t+\sigma_\theta\dd W,\qquad U(0)=U_0\in H .
\end{equation}


\subsection{Well-posedness and the path solution map}

We recall the well-posedness and ergodic properties established in~\cite{FGRT15}.
The pathwise regularity will also be used to define the controlled dynamics.

\begin{proposition}\label{prop:wp}
Assume (A1)-(A3) and let $U_0\in H$.
Then \eqref{eq:abstract} has a unique $\mathcal F_t$-adapted mild solution $U$.
Almost surely, $U\in C([0,\infty);H)\cap L^2_{loc}([0,\infty);H^1)$.
Moreover,
\begin{enumerate}
\item[(i)] for any $t_0>0$ and $m\ge0$, almost surely $U\in C([t_0,\infty);H^m)$;
\item[(ii)] for each noise realization, $U_0\mapsto U(t,U_0)$ is Fr\'echet differentiable on $H$.
For each $U_0$, the map $W\mapsto U(t,W)$ from $C([0,t];\R^4)$ to $H$ is Fr\'echet differentiable;
\item[(iii)] $(U_t,\Prob_U)_{U\in H}$ is a Feller Markov family on $H$, and it has a unique invariant measure, denoted by $\mu_*$.
\end{enumerate}
\end{proposition}

For $T>0$, set $|\zeta|_\infty:=\sup_{t\in[0,T]}|\zeta(t)|$ and
$C_0([0,T];\R^4):=\{\zeta\in C([0,T];\R^4):\zeta(0)=0\}$.
Given $\zeta\in C_0([0,T];\R^4)$, put $\xi(t):=\sigma_\theta\zeta(t)$.
This function has only a temperature component and is smooth in space.
Consider the shifted equation
\begin{equation}\label{eq:shifted}
\partial_t\bar U+A\bar U+B(\bar U+\xi)=G(\bar U+\xi)-A\xi,\qquad \bar U(0)=U_0 .
\end{equation}
Define the path solution map by $\Phi_t(U_0,\zeta):=\bar U(t)+\xi(t)$.
The next lemma (the proof is given in Section~\ref{sec:apriori}) gives its continuity and concatenation properties, which allow us to combine control paths in the irreducibility proof.

\begin{lemma}\label{lem:phi}
Assume (A1)-(A3). Then the following assertions hold:
\begin{enumerate}
\item[(i)] Let $U_0\in H$ and $\zeta\in C_0([0,T];\R^4)$.
Equation \eqref{eq:shifted} has a unique solution $\bar U\in C([0,T];H)\cap L^2(0,T;H^1)$.
For all $t_0>0$ and $m\ge0$, it satisfies $\bar U\in C([t_0,T];H^m)$.
Moreover, almost surely $U_t=\Phi_t(U_0,W|_{[0,t]})$ for all $t\in[0,T]$ simultaneously.
\item[(ii)] For any $t\in(0,T]$ and $R,\varrho>0$, there exists $C_*=C_*(t,R,\varrho)$ such that
\[
\norm{\Phi_t(U_0^1,\zeta^1)-\Phi_t(U_0^2,\zeta^2)}\le C_*\big(\norm{U_0^1-U_0^2}+|\zeta^1-\zeta^2|_\infty\big)
\]
whenever $\norm{U_0^i}\le R$ and $|\zeta^i|_\infty\le\varrho$ for $i=1,2$.
The constants can be chosen so that $C_*(s,R,\varrho)\le C_*(t,R,\varrho)$ for $s\le t$.
In particular, $\Phi_t:H\times C_0([0,t];\R^4)\to H$ is jointly continuous.
\item[(iii)] For $\zeta_1\in C_0([0,s];\R^4)$ and $\zeta_2\in C_0([0,t];\R^4)$, let
\[
(\zeta_1\oplus\zeta_2)(r):=\begin{cases}\zeta_1(r), & r\in[0,s],\\ \zeta_1(s)+\zeta_2(r-s), & r\in[s,s+t].\end{cases}
\]
Then $\Phi_{s+t}(U_0,\zeta_1\oplus\zeta_2)=\Phi_t(\Phi_s(U_0,\zeta_1),\zeta_2)$.
\end{enumerate}
The deterministic assertions in (i)-(iii) require only (A1) and (A2).
\end{lemma}

%

\subsection{Feynman-Kac semigroups and weight functions}

Let $\Prob_U$ be the solution law with initial state $U$, and let $\E_U$ denote its expectation.
Write $\Pt\psi(U):=\E_U\psi(U_t)$.
For $\nu\in\cP(H)$, set $\Prob_\nu:=\int_H\Prob_U\,\nu(\dd U)$ and denote the corresponding expectation by $\E_\nu$.
For a function $f$ and measure $\mu$, write $\ip f\mu:=\int_Hf(U)\,\mu(\dd U)$.
Let $B_H(\hat U,r):=\{U\in H:\norm{U-\hat U}\le r\}$ and $B_H(r):=B_H(0,r)$.
For $R>0$, define $X_R:=\{U\in H^1:\norm U_{H^1}\le R\}$.
We equip $X_R$ with the topology induced by $H$, in which it is compact.
For $\gamma,M>0$, define the class of initial measures
\begin{equation}\label{eq:Lambda}
\Lam(\gamma,M):=\Big\{\nu\in\cP(H):\ \int_He^{\gamma\norm U^2}\nu(\dd U)\le M\Big\}.
\end{equation}

For a Polish space $X$, let $C_b(X)$ denote the bounded continuous functions.
Let $\cM_+(X)$ denote the finite non-negative Borel measures.
We equip $\cP(X)$ and $\cM_+(X)$ with the weak topology, defined by convergence against $C_b(X)$.
In this topology, $\cP(X)$ is Polish.
Given a measurable weight $\fw:X\to[1,\infty)$, set
\[
\norm f_{L^\infty_\fw}:=\sup_{u\in X}\frac{|f(u)|}{\fw(u)}.
\]
The spaces $C_\fw(X)$ and $L^\infty_\fw(X)$ consist of continuous and measurable functions, respectively, with finite weighted norm.
Write $\cP_\fw(X):=\{\mu\in\cP(X):\ip{\fw}{\mu}<\infty\}$.
Define
\[
C_b^1(H):=\{V\in C_b(H):V\ \text{is continuously Fr\'echet differentiable},\ \norm{\nabla V}_\infty<\infty\}.
\] 
Throughout the paper, we use the following weight functions and constants:
\begin{equation}\label{eq:weights}
\fm_\varkappa(U):=e^{\varkappa\norm U^2},\quad \fw_m(U):=1+\norm U^{2m},\quad
\eta_1:=\frac{\kappa}{8\fB_0},\quad \gamma_0:=(1-e^{-\kappa/2})\,\eta_1 .
\end{equation}
For $V\in C_b(H)$, we define the Feynman-Kac semigroup
\begin{equation}\label{eq:FK}
\PV t f(U):=\E_U\big[\Xi_t^Vf(U_t)\big],\qquad \Xi_t^V:=\exp\Big(\int_0^tV(U_s)\dd s\Big),
\end{equation}
and write $\PV t{}^*:\cM_+(H)\to\cM_+(H)$ for its dual.
Its transition kernel is $\PV t(U,O):=(\PV t{}^*\delta_U)(O)$ for $O\in\cB(H)$.
The bound $\Xi_t^V\ge e^{-t\norm V_\infty}$ gives $\PV t(U,\cdot)\ge e^{-t\norm V_\infty}\Pt(U,\cdot)$.
A measure $\mu$ satisfying $\PV t{}^*\mu=\lambda^t\mu$ for all $t\ge0$ is an eigenmeasure.
A function $h$ satisfying $\PV th=\lambda^th$ for all $t\ge0$ is an eigenfunction.
We refer to these identities as eigenrelations.

\subsection{Main results}

Our first result is a multiplicative ergodic theorem for the Feynman-Kac semigroups.

\begin{theorem}\label{thm:B}
Assume (A1)-(A3) and let $V\in C_b^1(H)$.
Let $(m,R_0)$ be the constants associated with $V$ in Proposition~\ref{prop:growth1}.
Let $\eta_1$ be as in \eqref{eq:weights}.
Then,
\begin{enumerate}
\item[(i)] There exist $\lambda_V>0$, $h_V\in C_{\fw_m}(H)$ and $\mu_V\in\bigcap_{\varkappa\in(0,\eta_1)}\cP_{\fm_\varkappa}(H)$ such that for all $t\ge0$
\[
\PV th_V=\lambda_V^th_V,\qquad \PV t{}^*\mu_V=\lambda_V^t\mu_V,
\]
with $\ip{h_V}{\mu_V}=1$.
The function $h_V$ is strictly positive, and $\supp\mu_V=H$.
For each $\varkappa\in(0,\eta_1)$, $\mu_V$ is the unique eigenmeasure in $\cP_{\fm_\varkappa}(H)$.
Among eigenfunctions $h\in C_{\fw_m}(H)$ satisfying $\ip h{\mu_V}=1$, the function $h_V$ is unique.
\item[(ii)] For any $\varkappa\in(0,\eta_1)$, $\psi\in C_{\fm_\varkappa}(H)$ and $R>0$, we have
\[
\lambda_V^{-t}\PV t\psi\to\ip\psi{\mu_V}h_V\quad\text{in }C_b(B_H(R))\cap L^1(H,\mu_V)
\]
as $t\to\infty$.
For each $\nu\in\cP_{\fw_m}(H)$, we also have
\[
\lambda_V^{-t}\PV t{}^*\nu\to\ip{h_V}\nu\mu_V
\]
weakly in $\cM_+(H)$ as $t\to\infty$.
\item[(iii)] For any $\varkappa\in(0,\eta_1)$, $\varkappa'>\varkappa$, $M>0$ and $\psi\in C_{\fm_\varkappa}(H)$, it holds
\[
\lambda_V^{-t}\E_\nu\big[\Xi_t^V\psi(U_t)\big]\to\ip\psi{\mu_V}\ip{h_V}\nu
\]
as $t\to\infty$, uniformly with respect to $\nu\in\Lam(\varkappa',M)$.
\end{enumerate}
\end{theorem}

The second result concerns the occupation measures defined by
\[
\Lt:=\frac1t\int_0^t\delta_{U_s}\dd s,\qquad t>0 .
\]

\begin{theorem}\label{thm:A}
Assume (A1)-(A3). Then the following assertions hold:
\begin{enumerate}
\item[(i)] Let $\gamma>0$ and $M\ge1$.
The occupation measures satisfy a large deviation principle uniformly over $\nu\in\Lam(\gamma,M)$.
Its good rate function $I:\cP(H)\to[0,\infty]$ is independent of $\gamma$ and $M$.
For every weakly closed set $F\subset\cP(H)$ and weakly open set $G\subset\cP(H)$,
\begin{gather*}
\limsup_{t\to\infty}\frac1t\log\sup_{\nu\in\Lam(\gamma,M)}\Prob_\nu\{\Lt\in F\}\le-\inf_{\sigma\in F}I(\sigma),\\
\liminf_{t\to\infty}\frac1t\log\inf_{\nu\in\Lam(\gamma,M)}\Prob_\nu\{\Lt\in G\}\ge-\inf_{\sigma\in G}I(\sigma) .
\end{gather*}
\item[(ii)] $I(\sigma)=0$ if and only if $\sigma=\mu_*$, where $\mu_*$ is the unique invariant measure.
\item[(iii)] For every $\varepsilon\in(0,1]$, there exists $\sigma_\varepsilon\in\cP(H)$ such that
\[
\sigma_\varepsilon(B_H(0,\varepsilon))=1,
\qquad 0<I(\sigma_\varepsilon)<\infty.
\]
\item[(iv)] For $\sigma\in\cP(H)$,
\[
I(\sigma)=\sup_{V\in C_b(H)}\big(\ip V\sigma-Q(V)\big),\qquad
Q(V):=\lim_{t\to\infty}\frac1t\log\E_\nu\exp\Big(\int_0^tV(U_s)\dd s\Big),
\]
where the pressure limit is uniform over $\nu\in\Lam(\gamma,M)$ for each $\gamma>0$, $M\ge1$ and $V\in C_b(H)$.
Its value is independent of $\nu$, $\gamma$ and $M$, as is $I$.
Moreover, $Q$ is convex and $1$-Lipschitz on $C_b(H)$.
It satisfies $Q(V+c)=Q(V)+c$ for every $c\in\R$.
\end{enumerate}
\end{theorem}

The contraction principle gives the corresponding large deviation principle for time averages of bounded continuous observables.

\begin{corollary}\label{cor:level1}
Assume (A1)-(A3), and let $\gamma>0$, $M\ge1$ and $f\in C_b(H)$. Then the family
\[
\frac1t\int_0^tf(U_s)\dd s,\qquad t>0,
\]
satisfies the large deviation principle on $\R$, uniformly over $\nu\in\Lam(\gamma,M)$.
Its good rate function is
\[
I_f(x)=\inf\{I(\sigma):\sigma\in\cP(H),\ \ip f\sigma=x\},\qquad x\in\R .
\]
\end{corollary}

Note that Theorems~\ref{thm:B} and \ref{thm:A} and Corollary~\ref{cor:level1} are proved in Section~\ref{sec:proofs}.
\section{Well-posedness and a priori estimates}\label{sec:apriori}

The weighted energy inequality is the starting point for the estimates below. We first use it to justify the path solution map, and then derive the moment bounds needed for the Feynman-Kac semigroups.

\begin{lemma}\label{lem:energy}
Assume (A1)-(A3), let $U_0\in H$, and let $U=(w,\theta)$ be the solution given by Proposition~\ref{prop:wp}. Let
\[
M_t:=2\sum_{k\in\cZ,\,l\in\{0,1\}}\alpha_k^l\int_0^t\ip{U_s}{\sigma_k^l}_{\fc}\dd W_s^{k,l}.
\]
Then almost surely, for all $t\ge0$:
\begin{enumerate}
\item[(i)] the energy identity
\begin{equation}\label{energy-identity}
\norm{U_t}^2+2\int_0^tZ(U_s)\dd s=\norm{U_0}^2+2\fc g\int_0^t\ip{\partial_1\theta_s}{w_s}\dd s+\fB_0t+M_t
\end{equation}
holds;
\item[(ii)] the energy inequality
\begin{equation}\label{energy-inequality}
\norm{U_t}^2+\int_0^tZ(U_s)\dd s\le\norm{U_0}^2+\fB_0t+M_t
\end{equation}
holds,
where $M$ is a continuous local martingale, and
\begin{equation}\label{martingale}
\langle M\rangle_t=4\sum_{k\in\cZ,\,l\in\{0,1\}}(\alpha_k^l)^2\int_0^t\ip{U_s}{\sigma_k^l}_{\fc}^2\dd s\le4\fB_0\int_0^t\norm{U_s}^2\dd s\le\frac{4\fB_0}{\kappa}\int_0^tZ(U_s)\dd s;
\end{equation}
\item[(iii)] the deterministic solution $\Phi_t(U_0,0)$ with zero control satisfies $(i)$ and $(ii)$ with $\fB_0=0$ and $M\equiv0$. Consequently, $\norm{\Phi_t(U_0,0)}\le e^{-\kappa t/2}\norm{U_0}$.
\end{enumerate}
\end{lemma}

\begin{proof}
By the spatial regularity in Proposition~\ref{prop:wp}$(i)$ and the finite rank of the noise, It\^o's formula applies to $\norm{U_t}^2$ on $[t_0,T]$ for any $0<t_0<T$ and gives the energy identity on this interval.
Since $U\in C([0,\infty);H)\cap L^2_{loc}([0,\infty);H^1)$, letting $t_0\downarrow0$, we get \eqref{energy-identity}.
By the Young and Poincar\'e inequalities, we obtain
\begin{equation}\label{eq:buoy}
2\fc|g|\,|\ip{\partial_1\theta}{w}|\le\nu_2\norm{\nabla\theta}^2+\frac{\fc^2g^2}{\nu_2}\norm w^2=\nu_2\norm{\nabla\theta}^2+\fc\nu_1\norm w^2\le Z(U).
\end{equation}
Combining \eqref{eq:buoy} with \eqref{energy-identity} gives \eqref{energy-inequality}.
Independence of the Brownian coordinates gives the quadratic variation in \eqref{martingale}, whose bounds follow from Cauchy-Schwarz and $\kappa\norm U^2\le Z(U)$.
For zero control, the same computation with $\fB_0=0$ and $M=0$, followed by Gronwall's inequality, proves $(iii)$.
\end{proof}

We use the energy inequality to establish the properties of the path solution map.

\begin{proof}[Proof of Lemma~\ref{lem:phi}]
For $(i)$, since $\xi_w=0$, we have $B(\bar U+\xi)=B(\bar U)+B(\bar U,\xi)$, where the additional term is linear in $\bar U$ with smooth spatial coefficients.
Since the buoyancy term is controlled by \eqref{eq:buoy} and $G\xi-A\xi\in C([0,T];H^m)$, the Galerkin and regularity arguments of~\cite[Theorems~2.1.13, 2.1.18 and~2.1.19]{KS12} give the well-posedness of \eqref{eq:shifted} and the regularity at positive times.
The identity $U_t=\Phi_t(U_0,W|_{[0,t]})$ follows from the mild formulation and
\[
\int_0^te^{-(t-s)A}\sigma_\theta\dd W_s=\sigma_\theta W_t-\int_0^tAe^{-(t-s)A}\sigma_\theta W_s\dd s,
\]
where the integration by parts takes place on the finite-dimensional smooth noise range.

To prove $(ii)$, we write $U^i=(w^i,\theta^i):=\bar U^i+\xi^i$, $D:=\bar U^1-\bar U^2$ and $\varsigma:=\xi^1-\xi^2=\sigma_\theta(\zeta^1-\zeta^2)$.
Then
\[
\partial_tD+AD+B(D+\varsigma,U^1)+B(U^2,D+\varsigma)=G(D+\varsigma)-A\varsigma,
\]
where $B(\varsigma,U^1)=0$ because $\varsigma_w=0$, and $\ip{B(U^2,D)}D_\fc=0$.
Taking the weighted inner product with $D$ and using Ladyzhenskaya's inequality for $\ip{B(D,U^1)}D_\fc$, and \eqref{eq:buoy} and $\norm{\varsigma}_{W^{2,\infty}}\le C|\zeta^1-\zeta^2|_\infty$ for the remaining terms, we get
\[
\frac{\dd}{\dd t}\norm D^2\le C\big(1+\norm{U^1_t}_{H^1}^2\big)\norm D^2+C\big(1+\norm{U^2_t}^2\big)|\zeta^1-\zeta^2|_\infty^2 .
\]
The same energy computation for \eqref{eq:shifted} gives
\[
\norm{\bar U(t)}^2+\int_0^tZ(\bar U(s))\dd s\le\big(\norm{U_0}^2+Ct|\zeta|_\infty^2\big)e^{C(1+|\zeta|_\infty)t}.
\]
Since $\Phi_t(U_0^1,\zeta^1)-\Phi_t(U_0^2,\zeta^2)=D(t)+\varsigma(t)$, combining these two inequalities with Gronwall's inequality, we obtain $(ii)$ with a constant $C_*$ which depends only on $t,R,\varrho$ and is non-decreasing in $t$.

To prove $(iii)$, let $U_r:=\Phi_r(U_0,\zeta_1\oplus\zeta_2)$ and $\xi:=\sigma_\theta(\zeta_1\oplus\zeta_2)$.
Since the equation $\partial_r(U-\xi)=F(U)$ involves only the increments of the driving path, the functions $r\mapsto U_{s+r}$ and $r\mapsto\Phi_r(\Phi_s(U_0,\zeta_1),\zeta_2)$ solve the same problem on $[0,t]$ with the initial condition $U_s=\Phi_s(U_0,\zeta_1)$.
Hence the uniqueness in $(i)$ gives the required identity.
\end{proof}

We now record exponential estimates for the state, its time integral and discrete sums.
Their different forms will be used in the recurrence and Malliavin arguments.
\begin{lemma}\label{lem:expmom}
Assume (A1)-(A3) and let $U\in H$.
Set $q:=e^{-\kappa/2}$ and $c_q:=(1-q)^{-1}$.
Let $\eta_1$ and $\gamma_0$ be as in \eqref{eq:weights}.
Then,
\begin{enumerate}
\item[(i)] For $\eta\in(0,\eta_1]$ and $T>0$, 
\[
\E_U\exp\Big(\eta\norm{U_T}^2+\frac{\eta\kappa}{4}e^{-\kappa T/2}\int_0^T\norm{U_s}_{H^1}^2\dd s\Big)\le e^{1/4}\exp\big(\eta e^{-\kappa T/2}\norm U^2\big).
\]
\item[(ii)] For $\eta\in(0,\eta_1]$ and $T>0$,
\[
\E_U\exp\Big(\eta\norm{U_T}^2+\frac{\eta\kappa}{2}\int_0^T\norm{U_s}_{H^1}^2\dd s\Big)\le\exp\big(\eta\norm U^2+\eta\fB_0T\big).
\]
\item[(iii)] For $\eta\in(0,\gamma_0]$ and $N\in\N$,
\[
\E_U\exp\Big(\eta\sum_{k=0}^N\norm{U_k}^2\Big)\le\exp\big(c_q\eta\norm U^2+N/4\big).
\]
\item[(iv)] For $\eta\in(0,\eta_1]$ and $T>0$,
\[
\E_U\exp\Big(\frac\eta2\sup_{t\in[0,T]}\norm{U_t}^2\Big)\le2\exp\Big(\frac\eta2\norm U^2+\frac\eta2\fB_0T\Big).
\]
\end{enumerate}
\end{lemma}
\begin{proof}
Let $Y_t:=\norm{U_t}^2$ and $Z_t:=Z(U_t)$.
By \eqref{energy-identity}, \eqref{eq:buoy} and \eqref{martingale},
\begin{equation}\label{eq:Ysemimart}
Y_t=Y_0+\int_0^ta_s\dd s+M_t,\qquad a_s\le-Z_s+\fB_0,\qquad \dd\langle M\rangle_s\le\frac{4\fB_0}{\kappa}Z_s\dd s ,
\end{equation}
where $a_s:=-2Z_s+2\fc g\ip{\partial_1\theta_s}{w_s}+\fB_0$.
For $(i)$, set
\[
r(s):=\eta e^{-\kappa(T-s)/2},\qquad
\mathcal L_t:=r(t)Y_t+\frac14\int_0^tr(s)Z_s\dd s-\fB_0\int_0^tr(s)\dd s.
\]
Using \eqref{eq:Ysemimart}, $r'=\kappa r/2$ and $\kappa Y\le Z$, It\^o's formula shows that $e^{\mathcal L_t}$ is a non-negative local supermartingale.
Here $r\le\eta_1=\kappa/(8\fB_0)$ ensures the required sign of the drift.
Hence, by localization and Fatou's lemma, $\E e^{\mathcal L_T}\le e^{\mathcal L_0}$.
Since $r\ge\eta e^{-\kappa T/2}$, $Z\ge\kappa\norm\cdot_{H^1}^2$ and $\fB_0\int_0^Tr(s)\dd s\le1/4$, this implies $(i)$.
For $(ii)$, set
\[
\mathcal L_t:=\eta Y_t+\frac{\eta}{2}\int_0^tZ_s\dd s-\eta\fB_0t.
\]
By \eqref{eq:Ysemimart} and It\^o's formula, the drift of $e^{\mathcal L_t}$ is bounded above by
\[
e^{\mathcal L_t}\Big(-\frac{\eta}{2}+\frac{2\eta^2\fB_0}{\kappa}\Big)Z_t\dd t\le0,
\]
since $\eta\le\eta_1$. Thus $\E e^{\mathcal L_T}\le e^{\eta\norm U^2}$. Using $Z_s\ge\kappa\norm{U_s}_{H^1}^2$ yields $(ii)$.

To prove $(iii)$, let $\eta_N:=\eta$ and $\eta_{k-1}:=\eta+q\eta_k$ for $1\le k\le N$, so that $\eta_k\le c_q\eta\le\eta_1$.
Applying $(i)$ with $T=1$ successively and conditioning on $\mathcal F_{N-1},\dots,\mathcal F_0$, we get
$\E_U\exp\big(\eta\sum_{k=0}^N\norm{U_k}^2\big)\le e^{N/4}e^{\eta_0\norm U^2}$, which implies $(iii)$.

Finally, let $S:=e^{\mathcal L}$ with $\mathcal L$ as in $(ii)$.
By the maximal inequality for non-negative supermartingales, $\E\sup_{t\le T}S_t^{1/2}\le2S_0^{1/2}$.
Combining this with the inequality $\mathcal L_t\ge\eta(Y_t-\fB_0T)$ for $t\in[0,T]$, we get $(iv)$.
\end{proof}
The following consequences control exponential weights under the semigroup and under invariant measures.
\begin{corollary}\label{cor:expmom}
Assume (A1)-(A3), and let $\varkappa\in(0,\eta_1]$ and $t>0$. Then the following assertions hold:
\begin{enumerate}
\item[(i)] $\E_U\fm_\varkappa(U_t)\le e^{1/4}\fm_\varkappa(U)^{e^{-\kappa t/2}}\le e^{1/4}\fm_\varkappa(U)$, and for any $A>0$,
\[
\E_U\fm_\varkappa(U_t)\le e^{1/4}\Big[e^{-\varkappa(1-e^{-\kappa t/2})A}\fm_\varkappa(U)+e^{\varkappa A}\Big];
\]
\item[(ii)] for $V\in C_b(H)$ and $f\in L^\infty_{\fm_\varkappa}(H)$,
\[
|\PV tf(U)|\le e^{t\norm V_\infty+1/4}\norm f_{L^\infty_{\fm_\varkappa}}\fm_\varkappa(U);
\]
if $\varkappa<\eta_1$, then $\PV t$ maps $C_{\fm_\varkappa}(H)$ into itself; moreover, $\sup_{t\in[0,1]}\norm{\PV t\mathbf 1}_\infty\le e^{\norm V_\infty}$;
\item[(iii)] if $\varkappa\in(0,\eta_1)$, then there exists $C=C(\varkappa,\kappa,\fB_0)$ such that $\E_U[\norm{U_1}^2\fm_\varkappa(U_1)]\le C\fm_\varkappa(U)$;
\item[(iv)] $\mu_*\in\cP_{\fm_\varkappa}(H)$; consequently, there exists $M>0$ such that $\mu_*\in\Lam(\varkappa,M)$.
\end{enumerate}
\end{corollary}

\begin{proof}
Assertion $(i)$ follows from Lemma~\ref{lem:expmom}$(i)$, with the second bound obtained by distinguishing the cases $\norm U^2\ge A$ and $\norm U^2<A$.

The bound in $(ii)$ follows from $(i)$ and $\Xi_t^V\le e^{t\norm V_\infty}$.
For continuity, let $U^n\to U$ and couple the solutions with the same Brownian motion.
By Lemma~\ref{lem:phi}, the paths converge almost surely and uniformly on $[0,t]$, so the weighted terminal variables also converge almost surely.
Applying $(i)$ with $\varkappa'\in(\varkappa,\eta_1]$ gives uniform integrability, and Vitali's theorem therefore yields continuity of $\PV tf$.

For $(iii)$, choose $\varepsilon>0$ such that $\varkappa+\varepsilon\le\eta_1$ and $(\varkappa+\varepsilon)e^{-\kappa/2}\le\varkappa$. By $x\le\varepsilon^{-1}e^{\varepsilon x}$ and Lemma~\ref{lem:expmom}$(i)$ with $T=1$ we obtain
\[
\E_U\big[\norm{U_1}^2\fm_\varkappa(U_1)\big]\le\varepsilon^{-1}\E_U\fm_{\varkappa+\varepsilon}(U_1)\le\varepsilon^{-1}e^{1/4}\fm_\varkappa(U).
\]

For $(iv)$, choose $t,A>0$ in $(i)$ such that
\[
\Pt\fm_\varkappa\le\tfrac12\fm_\varkappa+e^{1/4+\varkappa A}.
\]
Let $\fm^R:=\fm_\varkappa\wedge e^{\varkappa R^2}$.
Then $\Pt\fm^R\le\tfrac12\fm^R+e^{1/4+\varkappa A}$ on $B_H(R)$ and $\Pt\fm^R\le\fm^R$ on $H\setminus B_H(R)$.
Integrating these inequalities against $\mu_*$ and using its invariance, we get
\[
\frac12\int_{B_H(R)}\fm_\varkappa\dd\mu_*
\le e^{1/4+\varkappa A}.
\]
Letting $R\to\infty$, we obtain $\ip{\fm_\varkappa}{\mu_*}\le2e^{1/4+\varkappa A}$.
The same estimate holds for any invariant probability measure.
\end{proof}
Polynomial moments provide the weights for the recurrence estimates.
Their decay rate can be increased by choosing the order of the moment.
\begin{lemma}\label{lem:poly}
Assume (A1)-(A3). Then for every integer $m\ge1$ there exists $C_m=C_m(m,\kappa,\fB_0)$ such that for all $U\in H$ and $t\ge0$,
\[
\E_U\norm{U_t}^{2m}\le e^{-m\kappa t/2}\norm U^{2m}+C_m .
\]
In particular,
\begin{enumerate}
\item[(i)] $\E_U\fw_m(U_t)\le e^{-m\kappa t/2}\fw_m(U)+C_m+1$, and hence $P_1\fw_m\le e^{-m\kappa/2}\fw_m+C_m+1$;
\item[(ii)] $\E[\norm{U_{k+1}}^{2m}\,|\,\mathcal F_k]\le e^{-m\kappa/2}\norm{U_k}^{2m}+C_m$ for $k\ge0$;
\item[(iii)] $\sup_{t\ge0}\E_U\norm{U_t}^{2m}\le\norm U^{2m}+C_m$.
\end{enumerate}
\end{lemma}

\begin{proof}
We derive from \eqref{eq:Ysemimart}, \eqref{martingale}, It\^o's formula and Young's inequality that
\[
\dd Y^m\le\Big[-\frac{m\kappa}{2}Y^m+C_m'\Big]\dd t+mY^{m-1}\dd M,
\]
where $C_m'$ depends only on $m,\kappa,\fB_0$.
By localization and Gronwall's inequality, this implies the required bound with $C_m=2C_m'/(m\kappa)$.
Assertions $(i)$-$(iii)$ follow from the definition of $\fw_m$ and the Markov property.
\end{proof}

We also use the following spatial regularity estimate from~\cite[Lemma~A.1$(iii)$, equation~(A.5)]{FGRT15}, with $s=1$.

\begin{lemma}\label{lem:H1mom}
Assume (A1)-(A3). Then for any $p\ge2$, $\eta>0$ and $T>0$ there exists $C=C(p,\eta,T,\nu_1,\nu_2,g,\fB_0)$ such that for all $U\in H$,
\[
\E_U\sup_{t\in[T/2,T]}\norm{U_t}_{H^1}^p\le C\exp(\eta\norm U^2).
\]
\end{lemma}
The moment bounds also give time continuity for polynomially growing test functions.
\begin{lemma}\label{lem:timecont}
Assume (A1)-(A3), and let $V\in C_b(H)$ and $m\ge1$. If $\psi\in C_{\fw_m}(H)$ and $U\in H$, then $t\mapsto\PV t\psi(U)$ is continuous on $[0,\infty)$.
\end{lemma}

\begin{proof}
The argument of~\cite[Lemma~3.3]{Ner19} applies with the moment estimate supplied by Lemma~\ref{lem:poly}.
\end{proof}
\section{Feynman-Kac semigroups: irreducibility, recurrence and the uniform Feller property}\label{sec:fk}
In this section, we verify the semigroup properties required for multiplicative ergodicity. We first establish uniform irreducibility, recurrence and growth bounds, and construct eigenmeasures. We then prove the uniform Feller property through a Malliavin gradient estimate.

\subsection{Uniform irreducibility}\label{sec:irr}\label{subsec:irr}
The sets $X_R$ are increasing and compact in $H$ by weak lower semicontinuity of the $H^1$ norm and the Rellich-Kondrachov theorem.
Their union $X_\infty:=\bigcup_{R>0}X_R=H^1$ is dense in $H$.
By Proposition~\ref{prop:wp}$(i)$, $P_t(U,\cdot)$ is concentrated on $X_\infty$ for every $U\in H$ and $t>0$.

Next, we recall full support of the Wiener measure and record a small-ball lower bound for long time intervals.
\begin{lemma}\label{lem:wiener}
Let $W$ be a standard Brownian motion in $\R^4$. Then the following assertions hold:
\begin{enumerate}
\item[(i)] for any $T>0$, $\zeta\in C_0([0,T];\R^4)$ and $\delta>0$,
\[
\Prob\Big(\sup_{s\in[0,T]}|W_s-\zeta(s)|<\delta\Big)>0;
\]
\item[(ii)] for any $\delta>0$ there exist $c_1,c_2>0$ such that for all $t\ge1$,
\[
\Prob\Big(\sup_{s\in[0,t]}|W_s|\le\delta\Big)\ge c_1e^{-c_2t}.
\]
\end{enumerate}
\end{lemma}

\begin{proof}
Assertion $(i)$ is the Wiener support property stated in~\cite[Remark~2.4.3]{KS12}.

To prove $(ii)$, fix $|x|\le\delta/4$ and note that $x+W_s=(1-s)x+(W_s+sx)$ implies
\[
\Big\{\sup_{s\le1}|W_s+sx|<\delta/8\Big\}
\subset
\Big\{\sup_{s\le1}|x+W_s|\le\delta,\ |x+W_1|\le\delta/4\Big\}.
\]
On $\{\sup_{s\le1}|W_s|<\delta/8\}$, we have $|x\cdot W_1|\le\delta^2/32$ and $|x|^2/2\le\delta^2/32$.
The Cameron-Martin formula with $h(s)=sx$ and the symmetry under $W\mapsto-W$ therefore give
\begin{align*}
\Prob\Big(\sup_{s\in[0,1]}|W_s+sx|<\frac\delta8\Big)&=e^{-|x|^2/2}\,\E\Big[\mathbf 1_{\{\sup_{s\in[0,1]}|W_s|<\delta/8\}}\,e^{-x\cdot W_1}\Big]\\
&\ge e^{-\delta^2/16}\,\Prob\Big(\sup_{s\in[0,1]}|W_s|<\frac\delta8\Big)=:p,
\end{align*}
where $p>0$ is independent of $x$ by $(i)$.
Since the path returns within distance $\delta/4$ of the origin at the end of each successful unit interval, iteration by the Markov property gives
\[
\Prob\Big(\sup_{s\in[0,t]}|W_s|\le\delta\Big)\ge p^{\lceil t\rceil}\ge p\,e^{-t\log(1/p)} .
\qedhere
\]
\end{proof}

Next, we will combine dissipation, control paths and the Wiener support property to obtain uniform irreducibility.
\begin{proposition}\label{prop:irr}
Assume (A1)-(A3). Then for any $r_*\ge1$, $R\ge1$ and $r>0$ there exist $l=l(r_*,r,R)>0$ and $p=p(r_*,r)>0$ such that
\[
P_l(U_0,B_H(\hat U,r))\ge p\qquad\text{for all }U_0\in X_R,\ \hat U\in X_{r_*}.
\]
If $V\in C_b(H)$ and $V\ge0$, then the same bound holds with $P_l$ replaced by $P_l^V$.
\end{proposition}
The dependence $l=l(r_*,r,R)$ and $p=p(r_*,r)$ is essential in applying the uniform irreducibility criterion.
\begin{proof}
Let $C_1$ be the constant in Lemma~\ref{lem:poly} for $m=1$, and set $d:=(2(1+C_1))^{1/2}$.
We first show that the process enters $B_H(0,d)$ with probability at least $1/2$.
For $U_0\in X_R$ and $t_R:=(4/\kappa)\log R$, Lemma~\ref{lem:poly} and Chebyshev's inequality give
\[
\Prob_{U_0}(\norm{U_{t_R}}>d)\le d^{-2}\big(e^{-\kappa t_R/2}R^2+C_1\big)\le d^{-2}(1+C_1)=\tfrac12 ,
\]
that is, $P_{t_R}(U_0,B_H(0,d))\ge1/2$.

It remains to reach each target from $B_H(0,d)$ in a time independent of $R$.
By compactness, we may cover $X_{r_*}$ by $n=n(r_*,r)$ balls $B_H(\hat U_i,r/4)$.
We apply~\cite[Theorem~5.18]{GHM18} with zero deterministic heat source.
For a driving path $q$ in the unit-amplitude formulation of that theorem, set $\zeta^{k,l}=q^{k,l}/\alpha_k^l$.
The product norm used there is equivalent to the weighted norm on $H$.
For $g<0$, the transformation $(w,\theta,\zeta)\mapsto(w,-\theta,-\zeta)$ gives the same system with buoyancy coefficient $-g>0$.

By approximate controllability (see~\cite[Theorem~5.18]{GHM18}), there are $\zeta_i\in C_0([0,1];\R^4)$ such that 
$$
\norm{\Phi_1(0,\zeta_i)-\hat U_i}<r/4, \quad  1\le i\le n.
$$
Applying Lemma~\ref{lem:phi}$(ii)$ with $T=R=1$ and $\varrho=\max_i|\zeta_i|_\infty+1$, we obtain $\varepsilon,\delta\in(0,1]$, common to all $i$, such that
\[
\norm{\Phi_1(v,\zeta)-\Phi_1(0,\zeta_i)}<r/4
\quad\text{if }\norm v\le\varepsilon,\quad |\zeta-\zeta_i|_\infty\le\delta.
\]
By Lemma~\ref{lem:energy}$(iii)$, choose $T_0>0$ with $e^{-\kappa T_0/2}d\le\varepsilon/2$ and set $l_2:=T_0+1$.
Lemma~\ref{lem:phi}$(ii)$, with $T=T_0$, $R=d$ and $\varrho=1$, then provides a sufficiently small $\varepsilon_0\in(0,1]$.
For $\norm v\le d$ and $|\zeta|_\infty\le\varepsilon_0$, we have
\[
\norm{\Phi_{T_0}(v,\zeta)}\le\norm{\Phi_{T_0}(v,0)}+\norm{\Phi_{T_0}(v,\zeta)-\Phi_{T_0}(v,0)}\le\tfrac\varepsilon2+\tfrac\varepsilon2=\varepsilon.
\]

Let $\hat U\in X_{r_*}$ and $v\in B_H(0,d)$, let $i$ be such that $\hat U\in B_H(\hat U_i,r/4)$, and let
\[
E_i:=\Big\{\sup_{s\in[0,T_0]}|W_s|\le\varepsilon_0\Big\}\cap\Big\{\sup_{s\in[0,1]}\big|(W_{T_0+s}-W_{T_0})-\zeta_i(s)\big|\le\delta\Big\}.
\]
By Lemma~\ref{lem:phi}$(i)$, $(iii)$, the solution issued from $v$ satisfies $U_{T_0}=\Phi_{T_0}(v,W|_{[0,T_0]})$ and $U_{l_2}=\Phi_1(U_{T_0},W_{T_0+\cdot}-W_{T_0})$ almost surely. Hence, on $E_i$, we have $\norm{U_{T_0}}\le\varepsilon$ and
\[
\norm{U_{l_2}-\hat U}\le\norm{U_{l_2}-\Phi_1(0,\zeta_i)}+\norm{\Phi_1(0,\zeta_i)-\hat U_i}+\norm{\hat U_i-\hat U}<\tfrac{3r}4.
\]
This implies that $P_{l_2}(v,B_H(\hat U,r))\ge\Prob(E_i)$.
Since the two events in $E_i$ are determined by the increments of $W$ on $[0,T_0]$ and on $[T_0,T_0+1]$, respectively, they are independent. Hence, by Lemma~\ref{lem:wiener}$(i)$,
\[
\Prob(E_i)\ge\Prob\Big(\sup_{s\in[0,T_0]}|W_s|\le\varepsilon_0\Big)\cdot\min_{1\le j\le n}\Prob\Big(\sup_{s\in[0,1]}|W_s-\zeta_j(s)|\le\delta\Big)=:p_2>0,
\]
where $l_2$ and $p_2$ depend only on $r_*$ and $r$.
Setting $l:=t_R+l_2$ and $p:=p_2/2$ and using the Chapman-Kolmogorov relation, we obtain
\[
P_l(U_0,B_H(\hat U,r))\ge\int_{B_H(0,d)}P_{t_R}(U_0,\dd v)\,P_{l_2}(v,B_H(\hat U,r))\ge\frac12p_2 .
\]
For $V\ge0$, the inequality $\Xi_l^V\ge1$ gives $P_l^V(U_0,\cdot)\ge P_l(U_0,\cdot)$.
\end{proof}
\begin{remark}\label{rem:irr-unif}
The abstract multiplicative ergodic theorem requires a lower bound independent of $R$.
For signed $V$, the estimate $\PV l\ge e^{-l\norm V_\infty}P_l$ need not preserve this independence, because $l$ depends on $R$.
We therefore set $c:=\inf_HV$ and use $\tilde V:=V-c\ge0$.
This leaves the eigenmeasure and eigenfunction unchanged.
The semigroups and eigenvalues are related by
\[
P_t^{\tilde V}=e^{-ct}P_t^V,\qquad
\lambda_{\tilde V}=e^{-c}\lambda_V,\qquad
\lambda_{\tilde V}^{t}=e^{-ct}\lambda_V^{t}.
\]
\end{remark}

\subsection{Recurrence, growth properties and eigenmeasures}\label{sec:growth}
In this subsection, we first estimate the return time to a compact set. The moment order will be chosen according to the required exponential rate.
\begin{proposition}\label{prop:rec}
Assume (A1)-(A3). For $a>0$, let
\[
m_0(a):=\Big\lceil\max\Big\{\frac{4\log2}\kappa,\ \frac{16a}\kappa\Big\}\Big\rceil .
\]
Let $m \ge m_0(a)$ be an integer.
There exist an integer $R_0=R_0(a,m)\ge1$, a constant $C=C(a,m)>0$ and an $\{\mathcal F_k\}$-stopping time $\tau_*$.
The stopping time takes values in $\{1,2,\dots\}\cup\{\infty\}$.
It is defined from $(U_k)_{k\in\N}$ by a formula independent of the initial state.
For every $U\in H$, we have $\tau_*<\infty$ and $U_{\tau_*}\in X_{R_0}$ almost surely under $\Prob_U$, and 
\[
\E_U\exp(a\tau_{R_0})\le\E_U\exp(a\tau_*)
\le C\fw_m(U),
\]
where $\tau_{R_0}
:=\inf\{t\ge0:\ U_t\in X_{R_0}\}$.
\end{proposition}

\begin{proof}
Let $q_m:=e^{-m\kappa/2}$ and $Y_k:=\norm{U_k}^{2m}$, and let $C_m$ be the constant in Lemma~\ref{lem:poly}.
Throughout the proof, we use the strong Markov property of the discrete-time process $(U_k)_{k\in\N}$.

\emph{Step 1: Return to a bounded set.}
Set $r^{2m}:=C_m/q_m$, $\sigma:=\inf\{k\in\N:\norm{U_k}\le r\}$ and $\beta_0:=m\kappa/4$.
On $\{\norm{U_k}>r\}$, Lemma~\ref{lem:poly}$(ii)$ gives $\E[Y_{k+1}\mid\mathcal F_k]\le2q_mY_k$.
Since $2q_me^{\beta_0}\le1$ by our choice of $m$, the process $e^{\beta_0(k\wedge\sigma)}Y_{k\wedge\sigma}$ is a supermartingale, which yields
$\Prob_U(\sigma>k)\le r^{-2m}e^{-\beta_0k}\norm U^{2m}$.
Since $\beta_0':=\beta_0/2\ge2a$, we obtain
\begin{equation}\label{eq:rec1}
\E_Ue^{\beta_0'\sigma}\le1+(e^{\beta_0'}-1)\sum_{j=0}^{\infty}e^{\beta_0'j}\Prob_U(\sigma>j)\le C\fw_m(U).
\end{equation}

From this ball, the process enters $X_{R_0}$ at time $1$ with probability tending to one as $R_0\to\infty$.
For $\norm v\le r$, Lemma~\ref{lem:H1mom} with $T=p=2$ and $\eta=1$, together with Chebyshev's inequality, gives
\[
\Prob_v(U_1\notin X_{R_0})\le CR_0^{-2}e^{r^2}=:p_0(R_0),
\]
here $C$ is independent of $v$ and $R_0$.

\emph{Step 2: Successive trials to enter a compact set.}
Let $\sigma_1:=\sigma$, $\sigma_{j+1}:=\inf\{k\ge\sigma_j+1:\norm{U_k}\le r\}$, $\ell_j:=\sigma_{j+1}-\sigma_j$, and
\[
J:=\inf\{j\ge1:\ U_{\sigma_j+1}\in X_{R_0}\},\qquad \tau_*:=\sigma_J+1 .
\]
The Markov property and \eqref{eq:rec1} give
$\E[e^{2a(\ell_j-1)}\mid\mathcal F_{\sigma_j+1}]\le C\fw_m(U_{\sigma_j+1})$.
Hence each $\sigma_j$ is finite almost surely.
The construction also shows that $\tau_*$ is a stopping time:
\[
\{\tau_*=k\}=\bigcup_{j\ge1}
\{\sigma_j=k-1,\ U_{\sigma_i+1}\notin X_{R_0}\text{ for }i<j,\ U_k\in X_{R_0}\}\in\mathcal F_k.
\]
Since $\norm{U_{\sigma_j}}\le r$, conditioning on $\mathcal F_{\sigma_j}$ allows us to apply the bound $p_0(R_0)$.
Together with Cauchy-Schwarz and Lemma~\ref{lem:poly} with order $2m$, this yields
\[
\E\big[e^{2a\ell_j}\mathbf 1_{\{U_{\sigma_j+1}\notin X_{R_0}\}}\,\big|\,\mathcal F_{\sigma_j}\big]\le e^{2a}C\,\E_{U_{\sigma_j}}\big[\mathbf 1_{\{U_1\notin X_{R_0}\}}\fw_m(U_1)\big]\le e^{2a}Cp_0^{1/2}\Theta_r=:\bar q,
\]
where $\Theta_r:=\big(2(1+q_{2m}r^{4m}+C_{2m})\big)^{1/2}$. Choose an integer $R_0$ such that $\bar q\le1/2$.

\emph{Step 3: Exponential moment of the entrance time.}
Let $Q_0:=1$ and
\[
Q_n:=\E\Big[\mathbf 1_{\{J>n\}}\prod_{j=1}^ne^{2a\ell_j}\Big],\qquad n\ge1 .
\]
Since $\{J>n-1\}$ and $\prod_{j=1}^{n-1}e^{2a\ell_j}$ are $\mathcal F_{\sigma_n}$-measurable, conditioning at $\sigma_n$ gives
\[
Q_n=\E\Big[\mathbf1_{\{J>n-1\}}\prod_{j=1}^{n-1}e^{2a\ell_j}\,
\E\big[e^{2a\ell_n}\mathbf1_{\{U_{\sigma_n+1}\notin X_{R_0}\}}\,\big|\,\mathcal F_{\sigma_n}\big]\Big]
\le\bar q Q_{n-1}\le2^{-n}.
\]
Thus $\Prob_U(J>n)\le Q_n\to0$.
Thus $\tau_*<\infty$ and $U_{\tau_*}\in X_{R_0}$ almost surely, and
\[
\E_Ue^{2a(\sigma_J-\sigma_1)}=\sum_{n\ge1}\E_U\Big[\mathbf 1_{\{J=n\}}\prod_{j=1}^{n-1}e^{2a\ell_j}\Big]\le\sum_{n\ge1}Q_{n-1}\le2 .
\]
Finally, by the Cauchy-Schwarz inequality, \eqref{eq:rec1} and $\fw_m\ge1$, we get
\[
\E_Ue^{a\tau_{R_0}}\le\E_Ue^{a\tau_*}\le e^a\big(\E_Ue^{2a\sigma_1}\big)^{1/2}\big(\E_Ue^{2a(\sigma_J-\sigma_1)}\big)^{1/2}\le\sqrt2e^a\big(C\fw_m(U)\big)^{1/2}\le C'\fw_m(U).
\qedhere
\]
\end{proof}

For $\psi:H\to\R$ and $R\ge1$, write $\norm\psi_R:=\sup_{U\in X_R}|\psi(U)|$.
The next three propositions give growth bounds and eigenmeasures for bounded continuous potentials.
This generality will be needed in the equilibrium-state argument.
We begin with polynomial weights.

\begin{proposition}\label{prop:growth1}
Assume (A1)-(A3).
Let $V\in C_b(H)$ and $\bar\gamma:=\Osc(V)=\sup V-\inf V$.
Choose an integer $m\ge m_0(\max\{\bar\gamma,1\})$.
Take $R_0$ and $\tau_*$ from Proposition~\ref{prop:rec} with $a=\max\{\bar\gamma,1\}$ and this $m$.
Write $\fw:=\fw_m$.
There exists $C=C(\bar\gamma,m,R_0)$ such that, for all $t\ge0$,
\begin{enumerate}
\item[(i)] $\PV t\mathbf 1(U)\le C\fw(U)\norm{\PV t\mathbf 1}_{R_0}$, $U\in H$;
\item[(ii)] $\norm{\PV t\fw}_{L^\infty_\fw}\le C\norm{\PV t\mathbf 1}_{R_0}$;
\item[(iii)] $\sup_{t\in[0,1]}\norm{\PV t\mathbf 1}_\infty\le e^{\norm V_\infty}$.
\end{enumerate}
These bounds hold for every such $m$, with its corresponding $R_0$.
We will use this freedom in Subsection~\ref{subsec:equil}.
\end{proposition}

\begin{proof}
Assertion $(iii)$ follows from Corollary~\ref{cor:expmom}$(ii)$.
Since the shift $V\mapsto V-\inf V$ multiplies both sides of $(i)$ and $(ii)$ by the same factor, we may assume $0\le V\le\bar\gamma$.
In particular, $\PV t\mathbf1\ge1$ and $\PV s\mathbf1\le\PV t\mathbf1$ for $s\le t$.

For $(i)$, Proposition~\ref{prop:rec} gives $\E_Ue^{\bar\gamma\tau_*}\le C_1\fw(U)$ for a constant $C_1>0$.
On $\{\tau_*>t\}$, the weight is at most $e^{\bar\gamma\tau_*}$.
On $\{\tau_*=k\}\in\mathcal F_k$, we have $U_k\in X_{R_0}$ and $\Xi_k^V\le e^{\bar\gamma k}$.
The Markov property and monotonicity give
\[
\E_U[\Xi_t^V;\tau_*\le t]=\sum_{k\in\N,\ k\le t}\E_U\big[\mathbf 1_{\{\tau_*=k\}}\Xi_k^V(\PV{t-k}\mathbf 1)(U_k)\big]\le\norm{\PV t\mathbf 1}_{R_0}\E_Ue^{\bar\gamma\tau_*}.
\]
Adding the two bounds and using $\norm{\PV t\mathbf 1}_{R_0}\ge1$, we obtain $(i)$ with $C=2C_1$.

We turn next to $(ii)$. By Lemma~\ref{lem:poly}$(i)$, $P_s\fw\le(C_m+2)\fw$, and hence
\[
\norm{\PV t\fw}_{L^\infty_\fw}\le e^{\bar\gamma}(C_m+2)\norm{\PV{[t]}\fw}_{L^\infty_\fw},\qquad \norm{\PV{[t]}\mathbf 1}_{R_0}\le\norm{\PV t\mathbf 1}_{R_0},
\]
so it suffices to consider integer times.
Set $q:=e^{\bar\gamma-m\kappa/2}<1$, using $m\ge16\bar\gamma/\kappa$.
The Markov property and Lemma~\ref{lem:poly}$(i)$ give
$\PV k\fw\le q\PV{k-1}\fw+e^{\bar\gamma}(C_m+1)\PV{k-1}\mathbf1$.
Iterating and applying $(i)$, we obtain
\begin{align*}
\PV k\fw(U)&\le q^k\fw(U)+e^{\bar\gamma}(C_m+1)\sum_{j=0}^{k-1}q^j\PV{k-1-j}\mathbf 1(U)\\
&\le\Big[1+\frac{2C_1e^{\bar\gamma}(C_m+1)}{1-q}\Big]\fw(U)\norm{\PV k\mathbf 1}_{R_0}.
\end{align*}
\end{proof}
The following estimates concern exponential weights and the $H^1$ energy.
\begin{proposition}\label{prop:growth23}
Assume (A1)-(A3).
Let $V\in C_b(H)$, $\bar\gamma:=\Osc(V)$ and $(m,R_0)$ be as in Proposition~\ref{prop:growth1}.
For $\varkappa\in(0,\eta_1)$, set $\fm:=\fm_\varkappa$.
Define $\mathsf F(U):=\norm U_{H^1}^2$ on $H^1$ and $\mathsf F(U):=+\infty$ otherwise.
There exists $C=C(\bar\gamma,\varkappa,m,R_0)$ such that
\begin{enumerate}
\item[(i)] $\norm{\PV t\fm}_{L^\infty_\fm}\le C\norm{\PV t\mathbf 1}_{R_0}$ for all $t\ge0$;
\item[(ii)] $\norm{\PV t\mathsf F}_{L^\infty_\fm}\le C\norm{\PV t\mathbf 1}_{R_0}$ for all $t\ge1$.
\end{enumerate}
\end{proposition}

\begin{proof}
As in Proposition~\ref{prop:growth1}, we may assume $0\le V\le\bar\gamma$.
Corollary~\ref{cor:expmom}$(i)$ then reduces $(i)$ to integer times.
Split $\PV k\fm$ according to $\norm{U_k}^2\le A$ and $\norm{U_k}^2>A$.
Proposition~\ref{prop:growth1}$(i)$ bounds the first part by $C_1e^{\varkappa A}\fm\,\norm{\PV k\mathbf1}_{R_0}$.
Here we use $\fw_m\le C_{m,\varkappa}\fm$, and $C_1$ is independent of $A$ and $k$.
Conditioning at time $k-1$ and applying Corollary~\ref{cor:expmom}$(iii)$, we bound the second part by $A^{-1}e^{\bar\gamma}C_2\PV{k-1}\fm$, where $C_2$ is independent of $A,k$.
Choosing $A$ so that $q:=A^{-1}e^{\bar\gamma}C_2\le1/2$, we obtain
\[
\norm{\PV k\fm}_{L^\infty_\fm}\le C_1e^{\varkappa A}\norm{\PV k\mathbf 1}_{R_0}+q\norm{\PV{k-1}\fm}_{L^\infty_\fm}.
\]
Iteration, monotonicity and $\norm{\PV k\mathbf1}_{R_0}\ge1$ prove $(i)$, exactly as in Proposition~\ref{prop:growth1}.

For $(ii)$, Lemma~\ref{lem:H1mom} with $T=p=2$ and $\eta=\varkappa$ gives $\PV1\mathsf F\le e^{\bar\gamma}C\fm$.
For $t\ge1$, the semigroup property and $(i)$ give
\[
\PV t\mathsf F=\PV{t-1}(\PV 1\mathsf F)\le e^{\bar\gamma}C\PV{t-1}\fm\le C\norm{\PV t\mathbf 1}_{R_0}\fm .
\qedhere
\]
\end{proof}
We next construct eigenmeasures and establish their moment and tail bounds.
These are the measure theoretic inputs to the multiplicative ergodic theorem.
\begin{proposition}\label{prop:eigen}
Assume (A1)-(A3), and let $V\in C_b(H)$ and $t>0$. Then
\begin{enumerate}
\item[(i)] there exist $\mu_{t,V}\in\cP(H)$ and $\lambda_{t,V}\in[e^{-t\norm V_\infty},e^{t\norm V_\infty}]$ such that $\PV t{}^*\mu_{t,V}=\lambda_{t,V}\mu_{t,V}$;
\item[(ii)] for any $\varkappa\in(0,\eta_1)$ and $n\ge1$,
\[
\int_H\big(\norm U_{H^1}^n+\fm_\varkappa(U)\big)\,\mu_{t,V}(\dd U)<\infty ;
\]
\item[(iii)] for any $m\ge1$,
\[
\norm{\PV t\fw_m}_R\int_{H\setminus X_R}\fw_m(U)\,\mu_{t,V}(\dd U)\to0\qquad\text{as }R\to\infty .
\]
\end{enumerate}
Moreover, $(ii)$ and $(iii)$ hold for every eigenmeasure belonging to $\cP_{\fm_\varkappa}(H)$ for some $\varkappa\in(0,\eta_1)$.
\end{proposition}

\begin{proof}
Fix $\varkappa\in(0,\eta_1)$ and write $q_t:=e^{-\kappa t/2}$.
For $A\ge1$, define
\[
\mathscr G(\chi):=\frac{\PV t{}^*\chi}{\ip{\PV t\mathbf1}\chi},
\qquad
D_A:=\{\chi\in\cP(H):\ip{\fm_\varkappa}\chi\le A\}.
\]
The set $D_A$ is convex, weakly closed and contains $\delta_0$.
By Corollary~\ref{cor:expmom}$(ii)$, $\PV t$ preserves $C_b(H)$.
Since the denominator in $\mathscr G$ is at least $e^{-t\norm V_\infty}$, it follows that $\mathscr G$ is weakly continuous.

For $\chi\in D_A$, Corollary~\ref{cor:expmom}$(i)$ and Jensen's inequality give
\[
\ip{\fm_\varkappa}{\mathscr G(\chi)}\le e^{t\Osc(V)+1/4}
\Big(\ip{\fm_\varkappa}\chi\Big)^{q_t}
\le e^{t\Osc(V)+1/4}A^{q_t}\le A,
\]
provided $A\ge\exp((t\Osc(V)+1/4)/(1-q_t))$.
Fix such an $A$, so that $\mathscr G(D_A)\subset D_A$.
Lemma~\ref{lem:H1mom}, with $T=2t$, $p=2$ and $\eta=\varkappa$, also gives
\[
\int_H\norm U_{H^1}^2\,\mathscr G(\chi)(\dd U)\le e^{t\Osc(V)}\int_H\E_U\norm{U_t}_{H^1}^2\,\chi(\dd U)\le e^{t\Osc(V)}CA=:C_1.
\]
The constant $C_1$ is independent of $\chi\in D_A$.
Set
\[
\mathscr K:=\{\chi\in D_A:\chi(H\setminus X_R)\le C_1R^{-2}\text{ for all }R\ge1\}.
\]
The preceding estimates show that $\mathscr K$ is non-empty, weakly closed and convex, with $\mathscr G(\mathscr K)\subset\mathscr K$.
Moreover, the tail bound implies tightness, so Prokhorov's theorem gives compactness.
By viewing finite signed measures in the Hausdorff locally convex topology generated by $\chi\mapsto|\ip f\chi|$, $f\in C_b(H)$, 
its restriction to $\cP(H)$ is the weak topology.
The Tychonoff fixed point theorem~\cite{Tyc35} therefore gives $\mu\in\mathscr K$ with $\mathscr G(\mu)=\mu$.
Setting $\lambda:=\ip{\PV t\mathbf1}\mu$, we obtain $\PV t{}^*\mu=\lambda\mu$, with the stated bounds on $\lambda$ and $\ip{\fm_\varkappa}\mu\le A$.

For $(ii)$, suppose $\ip{\fm_{\varkappa_0}}\mu<\infty$ and choose $\varkappa'\in(0,\eta_1)$ with $\varkappa'q_t\le\varkappa_0$.
The eigenrelation and Corollary~\ref{cor:expmom}$(i)$ then give
\[
\ip{\fm_{\varkappa'}\wedge N}\mu=\lambda^{-1}\ip{\PV t(\fm_{\varkappa'}\wedge N)}\mu\le\lambda^{-1}e^{t\norm V_\infty+1/4}\ip{\fm_{\varkappa_0}}\mu .
\]
Let $N\to\infty$ and iterate from $\varkappa_0=\varkappa$ to cover every $\varkappa'\in(0,\eta_1)$.
For the spatial moments, we combine the same truncated eigenrelation with Lemma~\ref{lem:H1mom}, taking $T=2t$, $p=\max\{n,2\}$ and $\eta\in(0,\eta_1)$, to obtain
\[
\ip{\norm\cdot_{H^1}^n\wedge N}\mu\le\lambda^{-1}e^{t\norm V_\infty}\big(1+C\ip{\fm_\eta}\mu\big)<\infty ,
\]
here $C$ is independent of $N$.

For $(iii)$, Lemma~\ref{lem:poly}$(i)$ gives $\norm{\PV t\fw_m}_R\le e^{t\norm V_\infty}(2+R^{2m}+C_m)$.
For $n>4m$, Cauchy-Schwarz, Chebyshev's inequality and the moments in $(ii)$ yield
\[
\int_{H\setminus X_R}\fw_m\dd\mu\le\ip{\fw_m^2}\mu^{1/2}\,\mu\big(\norm U_{H^1}>R\big)^{1/2}\le CR^{-n/2}.
\]
The product in $(iii)$ is therefore bounded by $C(1+R^{2m})R^{-n/2}$, which tends to zero.
The same argument applies to every eigenmeasure with one finite exponential moment.
\end{proof}


\subsection{Malliavin gradient estimate and the uniform Feller property}\label{sec:malliavin}\label{subsec:grad}
We now turn to the uniform Feller property.
The argument uses a gradient estimate for $\PV t$, followed by the mean value inequality.
We decompose the linearized flow into a Malliavin control term and a decaying residual.
Integration by parts handles the control term.
We begin with the linearized operators and their conditional Malliavin estimates.

Let $D$ denote the Malliavin derivative with respect to $W$, and let $\delta$ denote the Skorokhod integral.
Write $D_s^j$ for the $j$th component at time $s$.
For a separable Hilbert space $K$, let $\mathbb D^{1,p}(K)$ be the $K$-valued Malliavin-Sobolev space.
We abbreviate $\mathbb D^{1,p}(\R)$ to $\mathbb D^{1,p}$, following~\cite[Sections~1.2-1.3]{Nua06}.
All identities and estimates involving $D_s$ hold for almost every $s$, almost surely.

For $0\le s\le t$ and $\xi\in H$ let $\cJ_{s,t}\xi:=\phi(t)$, where $\phi$ is the solution of the linearized equation
\begin{equation}\label{eq:linearized}
\partial_r\phi+A\phi+\nabla B(U_r)\phi=G\phi,\quad r\in[s,t],\qquad \phi(s)=\xi,
\end{equation}
and $\nabla B(U)\phi:=B(U,\phi)+B(\phi,U)$. By Proposition~\ref{prop:wp}, $\cJ_{0,t}=\nabla_{U_0}U(t)$. By \cite[(3.3)-(3.4)]{FGRT15} we obtain
\[
D_s^jU(T)=\cJ_{s,T}\sigma_\theta\mathbf e_j,\qquad j=1,\dots,4,\ s\le T.
\]
We define the control operator $\cAop_{s,t}:L^2([s,t];\R^4)\to H$ and the Malliavin matrix $\mathfrak M_{s,t}$ by
\[
\cAop_{s,t}v:=\int_s^t\cJ_{r,t}\sigma_\theta v(r)\dd r,\qquad \mathfrak M_{s,t}:=\cAop_{s,t}\cAop_{s,t}^*,
\]
where $(\cAop_{s,t}^*\xi)(r)=\sigma_\theta^*\cK_{r,t}\xi$ for $r\in[s,t]$.
Here $\cK_{r,t}:=\cJ_{r,t}^*$ is the backward solution operator.
The matrix $\mathfrak M_{s,t}$ is non-negative and self-adjoint on $H$.
Operator Malliavin derivatives are understood through their action on fixed deterministic vectors.
Their operator norms are controlled by the derivative formulas below.
Taking suprema over countable dense subsets of the unit balls also ensures measurability.

Let $\cJ^{(2)}_{s,t}$ denote the second-order linearization.
For $s=0$, this is the second Fr\'echet derivative of $U_0\mapsto U_t$.
Pathwise, it is a bounded bilinear operator from $H\times H$ to $H$.
For such an operator $\mathcal B$, write
\[
|||\mathcal B|||:=\sup_{\norm\xi\le1,\ \norm{\xi'}\le1}\norm{\mathcal B(\xi,\xi')}.
\]
By~\cite[(A.22)]{FGRT15} we have
\[
D^j_\tau\cJ_{s,t}\xi=\begin{cases}\cJ^{(2)}_{\tau,t}(\sigma_\theta\mathbf e_j,\cJ_{s,\tau}\xi), & s\le\tau,\\[2pt] \cJ^{(2)}_{s,t}(\cJ_{\tau,s}\sigma_\theta\mathbf e_j,\xi), & s>\tau.\end{cases}
\]
For each integer $n\ge0$, define $W^{(n)}_r:=W_{n+r}-W_n$.
This process is independent of $\mathcal F_n$, and $U_{n+r}=\Phi_r(U_n,W^{(n)}|_{[0,r]})$.
Thus the operators on $[n,n+T]$ depend on $(U_n,W^{(n)})$.

We collect the linearization estimates needed on each time block.
The pathwise bounds are~\cite[(A.12), (A.18)-(A.21), (3.34)]{FGRT15}.
The moment estimates are the conditional versions of~\cite[(A.13)-(A.15), (A.23)-(A.25)]{FGRT15}.
The proof below addresses this conditional formulation.

\begin{lemma}\label{lem:linmom}
Assume (A1)-(A3), and let $T\in\{1,2\}$, $n\in\N$, $\eta>0$ and $p\ge1$. Then there exists $C=C(\eta,p,T)$ such that the following assertions hold.
\begin{enumerate}
\item[(i)] For $0\le s\le t$,
\[
\norm{\cJ_{s,t}}_{L(H)}\le\exp\Big(\eta\int_s^t\norm{U_r}_{H^1}^2\dd r+C_\eta(t-s)\Big),
\]
where $C_\eta$ depends on $\eta,\nu_1,\nu_2,g$ and is independent of $s,t$.
\item[(ii)] We have
\[
\E\Big[\sup_{n\le s<t\le n+T}\big(\norm{\cJ_{s,t}}^p+\norm{\cK_{s,t}}^p+|||\cJ^{(2)}_{s,t}|||^p\big)\,\Big|\,\mathcal F_n\Big]\le Ce^{\eta\norm{U_n}^2}.
\]
\item[(iii)] For $0\le s<t$ and $\beta>0$,
\begin{gather*}
\norm{\cAop_{s,t}}_{L(L^2,H)}\le C\Big(\int_s^t\norm{\cJ_{r,t}}^2\dd r\Big)^{1/2},\qquad \norm{(\mathfrak M_{s,t}+\beta\,\mathrm{Id})^{-1/2}}\le\beta^{-1/2},\\
\norm{\cAop^*_{s,t}(\mathfrak M_{s,t}+\beta\,\mathrm{Id})^{-1/2}}\le1,\qquad \norm{(\mathfrak M_{s,t}+\beta\,\mathrm{Id})^{-1/2}\cAop_{s,t}}\le1.
\end{gather*}
\item[(iv)] For $n\le s<t\le n+T$, the operators $\cJ_{s,t}$, $\cAop_{s,t}$ and $\cAop^*_{s,t}$ are Malliavin differentiable in the sense above.
For almost every $\tau\in[n,t]$ and each $j\in\{1,\dots,4\}$,
\[
\E\Big[\norm{D^j_\tau\cJ_{s,t}}^p+\norm{D^j_\tau\cAop_{s,t}}^p_{L(L^2,H)}+\norm{D^j_\tau\cAop^*_{s,t}}^p_{L(H,L^2)}\,\Big|\,\mathcal F_n\Big]\le Ce^{\eta\norm{U_n}^2}.
\]
\item[(v)] The operator $(\mathfrak M_{s,t}+\beta\,\mathrm{Id})^{-1}$ is Malliavin differentiable, and
\[
D^j_\tau(\mathfrak M+\beta\,\mathrm{Id})^{-1}=-(\mathfrak M+\beta\,\mathrm{Id})^{-1}\big[(D^j_\tau\cAop)\cAop^*+\cAop(D^j_\tau\cAop^*)\big](\mathfrak M+\beta\,\mathrm{Id})^{-1}.
\]
\end{enumerate}
\end{lemma}

\begin{proof}
Assertions $(i)$, $(iii)$ and $(v)$ follow from~\cite[(A.12), (A.18)-(A.22)]{FGRT15} and
\[
D(\mathfrak M+\beta\,\mathrm{Id})^{-1}
=-(\mathfrak M+\beta\,\mathrm{Id})^{-1}(D\mathfrak M)
 (\mathfrak M+\beta\,\mathrm{Id})^{-1}.
\]
For $(ii)$, write $W^{(n)}_r=W_{n+r}-W_n$. The Markov property of the solution and variation equations yields, for every measurable functional $\Psi$ appearing in $(ii)$,
\[
\E[\Psi(U_n,W^{(n)})\mid\mathcal F_n]
=\left.\E[\Psi(u,W)]\right|_{u=U_n}
\le C e^{\eta\norm{U_n}^2},
\]
by~\cite[(A.13)-(A.15)]{FGRT15}.

For $(iv)$, first fix deterministic input vectors and $p>1$. The solutions and their first and second variations admit jointly measurable versions. A Cameron-Martin perturbation supported on $(n,n+T)$ leaves $U_n$ unchanged. The Malliavin chain rule and the identities in~\cite[(A.22)]{FGRT15} therefore give, for $n<\tau<n+T$,
\begin{equation}\label{eq:Dshift}
D_\tau X(U_n,W^{(n)})=
\left.(D^+_{\tau-n}X)(u,W^{(n)})\right|_{u=U_n},
\end{equation}
where $D^+$ differentiates only the future increments and $X$ denotes one of $\cJ$, $\cAop$ or $\cAop^*$.
By independence of $W^{(n)}$ and $\mathcal F_n$, Fubini's theorem and the Cameron-Martin characterization~\cite[Section~2.1]{Nua06}, the identities for deterministic initial data extend to $U_n$.
Choose $\eta_0\in(0,\eta_1]$ sufficiently small. The estimates~\cite[(A.23)-(A.25)]{FGRT15} and Lemma~\ref{lem:expmom}$(i)$ give
\[
\E\big[\norm{X(U_n,W^{(n)})}^p+
\norm{D^+X(U_n,W^{(n)})}^p\big]
\le C\E e^{\eta_0\norm{U_n}^2}<\infty.
\]
Hence, for almost every $\tau\in[n,n+T]$,
\[
\E\big[\norm{D_\tau X(U_n,W^{(n)})}^p\,\big|\,\mathcal F_n\big]
=\left.\E\norm{D^+_{\tau-n}X(u,W)}^p\right|_{u=U_n}
\le C e^{\eta\norm{U_n}^2}.
\]
The same conditioning applies to $X$ itself. Applying the conditional identity to the operator norms in~\cite[(A.23)-(A.25)]{FGRT15} gives $(iv)$.
These norms have measurable versions obtained from countable dense sets of input vectors.
The case $p=1$ follows from conditional H\"older's inequality.
\end{proof}

\begin{remark}
In the conditional estimate, $D_\tau U_n=0$ for $\tau>n$. The derivatives on $(n,n+T)$ are therefore derivatives of the restarted flow with the initial state held fixed. This distinction is used in the estimate of the Skorokhod integral in Lemma~\ref{lem:residual}.
\end{remark}

To decompose the linearized flow, we use the alternating control construction of~\cite[Section~3]{FGRT15}.
Fix $U_0=U\in H$, a unit vector $\xi\in H$ and $\beta>0$.
Set $\rho_0:=\xi$. For $k\ge0$, define
\begin{align}
v(r)&:=\begin{cases}\big(\cAop^*_{2k,2k+1}(\mathfrak M_{2k,2k+1}+\beta\,\mathrm{Id})^{-1}\cJ_{2k,2k+1}\rho_{2k}\big)(r), & r\in(2k,2k+1],\\[2pt] 0, & r\in(2k+1,2k+2],\end{cases}\label{eq:control}\\
\rho(t)&:=\cJ_{0,t}\xi-\cAop_{0,t}v,\qquad \rho_n:=\rho(n).\label{eq:residual}
\end{align}
Set $v(0)=0$.
By~\cite[(3.7)]{FGRT15}, the residual satisfies
\[
\partial_t\rho+A\rho+\nabla B(U_t)\rho=G\rho-\sigma_\theta v,\qquad \rho(0)=\xi.
\]
The identity $\mathfrak M=\cAop\cAop^*$ gives the recursion
\[
\rho_{2k+1}=\mathcal R^\beta_{2k,2k+1}\cJ_{2k,2k+1}\rho_{2k},\qquad \rho_{2k+2}=\cJ_{2k+1,2k+2}\rho_{2k+1},
\]
where $\mathcal R^\beta_{s,t}:=\beta(\mathfrak M_{s,t}+\beta\,\mathrm{Id})^{-1}$ and $\norm{\mathcal R^\beta_{s,t}}\le1$.
The next block estimate allows the residual contraction to be chosen arbitrarily strong.

\begin{lemma}\label{lem:FGRT31}
Assume (A1)-(A3) and fix $\varpi,\delta_*>0$.
There exists $\beta=\beta(\varpi,\delta_*)>0$ such that, for every even $n\ge0$, the corresponding residual satisfies
\begin{equation}\label{eq:FGRT316}
\E\big(\norm{\rho_{n+2}}^8\,\big|\,\mathcal F_n\big)\le\delta_*\exp(\varpi\norm{U_n}^2)\norm{\rho_n}^8 .
\end{equation}
Moreover, $\rho(t)$ and $v(t)$ are $\mathcal F_{t^\sharp}$-measurable, where
\[
t^\sharp:=\begin{cases}\lceil t\rceil, & \lceil t\rceil\ \text{odd},\\ t, & \lceil t\rceil\ \text{even}.\end{cases}
\]
\end{lemma}

\begin{proof}
The construction above satisfies the assumptions of~\cite[Lemma~3.1]{FGRT15}, so the assertion follows from its equations~(3.16)-(3.18).
\end{proof}

By~\cite[Lemma~A.7]{FGRT15}, $\rho_{2k}\in\mathbb D^{1,p}(H)$ for every $p>1$.
The control satisfies $v\in\mathbb D^{1,2}(L^2([0,N];\R^4))$ for every $N>0$.
Thus $v\mathbf1_{[0,t]}\in\mathrm{Dom}\,\delta$ by~\cite[Proposition~1.3.1]{Nua06}.
Throughout, $\int_0^tv(s)\cdot\dd W_s$ denotes the Skorokhod integral $\delta(v\mathbf1_{[0,t]})$.

We now derive estimates with a prescribed decay rate for the residual and the cost of control.
\begin{lemma}\label{lem:residual}
Assume (A1)-(A3) and fix $\eta,\mathfrak a>0$.
There exist $\beta=\beta(\eta,\mathfrak a)>0$ and $C=C(\eta,\mathfrak a)>0$ with the following properties.
For every $U\in H$ and unit vector $\xi\in H$, construct $(v,\rho)$ by \eqref{eq:control}-\eqref{eq:residual} using this $\beta$.
Then,
\begin{enumerate}
\item[(i)] $\E_U\norm{\rho(t)}^4\le C\exp(\eta\norm U^2-\mathfrak a t)$ for $t\ge0$;
\item[(ii)] for every integer $n\ge0$ and $t\in[n,n+1]$,
\[
\E_U\Big|\int_n^tv(s)\cdot\dd W_s\Big|^2\le C\exp(\eta\norm U^2-\mathfrak a n);
\]
\item[(iii)] The control vanishes on $(2k+1,2k+2]$.
The integral $\int_{2k}^{2k+1}v(s)\cdot\dd W_s$ is $\mathcal F_{2k+1}$-measurable.
The residual $\rho_{2k}$ is $\mathcal F_{2k}$-measurable, and $D_s\rho_{2k}=0$ for $s>2k$.
\end{enumerate}
\end{lemma}

\begin{proof}
\emph{Weighted endpoint estimate}. Following~\cite[Section~3.4]{FGRT15}, set $c_q:=(1-e^{-\kappa/2})^{-1}$ and
\[
\mathfrak a_0:=2\mathfrak a,\quad \varpi:=\min\{\gamma_0,2\eta/c_q\},\quad \eta':=\varpi/2,\quad \delta_*:=e^{-4\mathfrak a_0-1/2},\quad \beta:=\beta(\varpi,\delta_*),
\]
where $\beta(\varpi,\delta_*)$ is given by Lemma~\ref{lem:FGRT31}.

Since a zero residual remains zero, \eqref{eq:FGRT316} yields
\[
\E\Big[\frac{\norm{\rho_{2k+2}}^8}{\norm{\rho_{2k}}^8}
 e^{-\varpi\norm{U_{2k}}^2}\,\Big|\,\mathcal F_{2k}\Big]\le\delta_*,
\]
where the ratio is defined to be zero when $\rho_{2k}=0$.
The factors preceding the $k$th block are $\mathcal F_{2k}$-measurable. Iterated conditioning thus gives
\[
\E\prod_{k=0}^{n-1}
\Big(\frac{\norm{\rho_{2k+2}}^8}{\norm{\rho_{2k}}^8}
 e^{-\varpi\norm{U_{2k}}^2}\Big)\le\delta_*^n.
\]
Since $2\eta'=\varpi$, Cauchy-Schwarz and Lemma~\ref{lem:expmom}$(iii)$ imply
\begin{align*}
\E\big[\norm{\rho_{2n}}^4e^{\eta'\norm{U_{2n}}^2}\big]
&\le\delta_*^{n/2}
\Bigg(\E\exp\Big(\varpi\sum_{k=0}^{2n}\norm{U_k}^2\Big)\Bigg)^{1/2}\\
&\le\delta_*^{n/2}e^{c_q\varpi\norm U^2/2+n/4}.
\end{align*}
Since $c_q\varpi/2\le\eta$ and $\delta_*^{1/2}e^{1/4}=e^{-2\mathfrak a_0}$, we obtain
\begin{equation}\label{eq:rho-even}
\E\big[\norm{\rho_{2n}}^4e^{\eta'\norm{U_{2n}}^2}\big]\le\exp(\eta\norm U^2-2n\mathfrak a_0),\qquad n\ge0.
\end{equation}

\emph{Residual estimate within a block}. Set $S_n:=\sup_{2n\le s<t\le2n+2}\norm{\cJ_{s,t}}$.
Lemma~\ref{lem:linmom}$(iii)$ gives $\norm v_{L^2([2n,2n+1])}\le\beta^{-1/2}S_n\norm{\rho_{2n}}$ and $\norm{\cAop_{2n,t}}\le CS_n$.
Inserting these bounds into the residual equation and using the vanishing control on $(2n+1,2n+2]$, we obtain
\[
\norm{\rho(t)}\le C(1+\beta^{-1/2})(1+S_n^3)\norm{\rho_{2n}},\qquad t\in[2n,2n+2].
\]
Conditioning on $\mathcal F_{2n}$, we combine Lemma~\ref{lem:linmom}$(ii)$ with \eqref{eq:rho-even}, taking $T=2$, $p=12$ and exponent $\eta'$, to obtain
\[
\E\norm{\rho(t)}^4\le C_\beta\E\big[\norm{\rho_{2n}}^4\E\big((1+S_n^3)^4\,\big|\,\mathcal F_{2n}\big)\big]\le C_\beta\E\big[\norm{\rho_{2n}}^4e^{\eta'\norm{U_{2n}}^2}\big]\le C_\beta e^{2\mathfrak a_0}e^{\eta\norm U^2-\mathfrak a_0t},
\]
where $C_\beta$ depends only on $\beta,\eta'$. Since $\mathfrak a_0\ge\mathfrak a$, $(i)$ follows.

\emph{Skorokhod integral of the control}. Since the control vanishes on odd blocks, consider $n=2k$ and $t\in[2k,2k+1]$.
We omit the operator subscripts $2k,2k+1$ in the following estimates.
Extend $u:=v\mathbf1_{[2k,t]}$ by zero to $[0,2k+1]$.
The generalized It\^o isometry~\cite[(1.54)]{Nua06} gives
\begin{align*}
\E|\delta(u)|^2
&=\E\int_0^{2k+1}|u(r)|^2\dd r+\sum_{i,j=1}^4\E\int_0^{2k+1}\int_0^{2k+1}
D_s^iu^j(r)\,D_r^ju^i(s)\dd s\dd r.
\end{align*}
Since $u(r)=0$ outside $[2k,t]$, its Malliavin derivatives vanish there as well.
Thus the product in the double integral vanishes unless both $r$ and $s$ lie in $[2k,t]$.
By Cauchy-Schwarz and symmetry in $(r,s,i,j)$, we obtain
\begin{align*}
\E\Big|\int_{2k}^tv(s)\cdot\dd W_s\Big|^2
\le\underbrace{\E\int_{2k}^{2k+1}|v(s)|^2\dd s}_{L_1}+\underbrace{\sum_{j=1}^4\int_{2k}^{2k+1}
\E\norm{D^j_sv}^2_{L^2([2k,2k+1];\R^4)}\dd s}_{L_2}.
\end{align*}
Since $\rho_{2k}$ is $\mathcal F_{2k}$-measurable, we have $D_s\rho_{2k}=0$ for $s>2k$ by~\cite[Corollary~1.2.1]{Nua06}. The resolvent bounds in Lemma~\ref{lem:linmom}$(iii)$ give
\[
\norm{v}_{L^2([2k,2k+1])}\le\beta^{-1/2}\norm{\cJ_{2k,2k+1}}\norm{\rho_{2k}}.
\]
By the product rule, Lemma~\ref{lem:linmom}(v) and~\cite[(3.33)-(3.35)]{FGRT15}, for almost every $s\in[2k,2k+1]$ we have
\[
\norm{D_s^jv}_{L^2}\le\Big(\beta^{-1/2}\norm{D_s^j\cJ}+\beta^{-1}\norm{D_s^j\cAop}\norm{\cJ}+2\beta^{-1}\norm{D_s^j\cAop^*}\norm{\cJ}\Big)\norm{\rho_{2k}}.
\]
Using the H\"older inequality and Lemma~\ref{lem:linmom}$(ii)$, $(iv)$ (with fourth moments and exponential parameters whose sum is $\eta'$), we obtain
\[
\E\Big[\norm v_{L^2}^2+\sum_{j=1}^4\int_{2k}^{2k+1}\norm{D_s^jv}_{L^2}^2\dd s\,\Big|\,\mathcal F_{2k}\Big]
\le C(\beta^{-1}+\beta^{-2})e^{\eta'\norm{U_{2k}}^2}\norm{\rho_{2k}}^2.
\]
Combining this estimate with \eqref{eq:rho-even}, Corollary~\ref{cor:expmom}$(i)$ and Cauchy-Schwarz, we deduce
\begin{align*}
L_1+L_2&\le C_\beta\E\big[e^{\eta'\norm{U_{2k}}^2}\norm{\rho_{2k}}^2\big]\\
&\le C_\beta\Big(\E[\norm{\rho_{2k}}^4e^{\eta'\norm{U_{2k}}^2}]\Big)^{1/2}
\Big(\E e^{\eta'\norm{U_{2k}}^2}\Big)^{1/2}\\
&\le C_\beta e^{\eta\norm U^2-k\mathfrak a_0}=C_\beta e^{\eta\norm U^2-\mathfrak a n}.
\end{align*}
This proves $(ii)$.

\emph{Measurability}. The block measurability and vanishing future derivative  in $(iii)$ follow from~\cite[(3.17)-(3.18), (A.30)]{FGRT15}.
For the Skorokhod integral, set $u:=v\mathbf1_{[2k,2k+1]}$.
Since $u$ is $\mathcal F_{2k+1}$-measurable and supported in $[2k,2k+1]$, the duality formula and~\cite[Proposition~1.2.8]{Nua06} give
\[
\E[\mathsf X\delta(u)]
=\E\big[\E(\mathsf X\mid\mathcal F_{2k+1})\delta(u)\big],
\qquad \mathsf X\in\mathbb D^{1,2}.
\]
Density in $L^2(\Omega)$ implies that $\delta(u)$ is $\mathcal F_{2k+1}$-measurable.
\end{proof}

\begin{remark}
The dependence $\beta=\beta(\eta,\mathfrak a)$ in Lemma~\ref{lem:residual} allows the decay exponent $\mathfrak a$ to be chosen after the potential $V$ is fixed. This is used in Proposition~\ref{prop:grad} to make both $e^{\bar\gamma-\mathfrak a/4}<1$ and $e^{2\bar\gamma-\mathfrak a}<1$.
\end{remark}

The gradient estimate now follows by combining the residual and control bounds of Lemma~\ref{lem:residual} with Malliavin integration by parts.

\begin{proposition}\label{prop:grad}
Assume (A1)-(A3) and let $V,\psi\in C_b^1(H)$.
Take $(m,R_0)$ from Proposition~\ref{prop:growth1} for this $V$.
For each $\eta>0$ and $\varepsilon_*\in(0,1)$, there exists $C=C(\eta,\varepsilon_*,\Osc(V),m,R_0)$ such that,
 for all $U,\xi\in H$ and $t\ge0$,
\[
|\nabla_\xi\PV t\psi(U)|\le Ce^{\eta\norm U^2}\norm{\PV t\mathbf 1}_{R_0}\big[(1+\norm{\nabla V}_\infty)\norm\psi_\infty+\varepsilon_*^t\norm{\nabla\psi}_\infty\big]\norm\xi .
\]
\end{proposition}

\begin{proof}
As in the proof of Proposition~\ref{prop:growth1}, replacing $V$ by $V-\inf V$, we can assume that $0\le V\le\bar\gamma:=\Osc(V)$.
Then $\PV t\mathbf1\ge1$ and $\PV s\mathbf1\le\PV t\mathbf1$ for $s\le t$.
By homogeneity, it suffices to consider the case $\norm\xi=1$.
Let
\[
\mathfrak a:=4\bar\gamma-4\log\varepsilon_*+1,
\qquad e^{\bar\gamma-\mathfrak a/4}=e^{-1/4}\varepsilon_*<\varepsilon_*.
\]
Use $\beta=\beta(\eta,\mathfrak a)$ from Lemma~\ref{lem:residual} to construct $(v,\rho)$ by \eqref{eq:control}-\eqref{eq:residual}, with $\rho_0=\xi$.
In particular, $\mathfrak a>4\bar\gamma$ and $2\bar\gamma-\mathfrak a<0$.

\emph{Step 1: Malliavin decomposition of the gradient.}
By Lemma~\ref{lem:linmom}$(i)$ with exponent $\eta_1\kappa/4$ and Lemma~\ref{lem:expmom}$(ii)$ with exponent $\eta_1$, the difference quotients of $U\mapsto\Xi_t^V\psi(U_t)$ are bounded in $L^2$, locally uniformly in $U$.
Hence, as in~\cite[(3.2)]{NPX23}, Proposition~\ref{prop:wp}$(ii)$ and Vitali's theorem give
\begin{equation}\label{eq:gradformula}
\nabla_\xi\PV t\psi(U)=\E_U\Big[\Xi_t^V\psi(U_t)\int_0^t\nabla V(U_s)\cJ_{0,s}\xi\dd s+\Xi_t^V\nabla\psi(U_t)\cJ_{0,t}\xi\Big].
\end{equation}

To decompose this expression, we follow~\cite[Section~3]{NPX23}.
By~\cite[(3.3)-(3.4), Lemma~A.7]{FGRT15}, we have $U_s\in\mathbb D^{1,2}(H)$ and $\cAop_{0,s}v=\ip{DU_s}v_{L^2([0,s];\R^4)}$, and this identity holds also for the non-adapted control $v$.
On each finite time interval, the moment bounds for $DU_s$ and the boundedness of $\nabla V$ allow us to differentiate $\int_0^tV(U_s)\dd s$ under the integral by closedness of the Malliavin derivative. Since this integral is bounded, the chain and product rules (see~\cite[Proposition~1.2.3]{Nua06}) give
\[
\ip{D(\Xi_t^V\psi(U_t))}{v}=\Xi_t^V\psi(U_t)\int_0^t\nabla V(U_s)\cAop_{0,s}v\dd s+\Xi_t^V\nabla\psi(U_t)\cAop_{0,t}v .
\]
Substituting $\cJ_{0,s}\xi=\cAop_{0,s}v+\rho(s)$ into \eqref{eq:gradformula}, we apply the duality formula~\cite[(1.42)]{Nua06}, since $v\mathbf1_{[0,t]}\in\mathrm{Dom}\,\delta$.
This gives $\nabla_\xi\PV t\psi(U)=I_1+I_2+I_3$, where
\begin{align*}
I_1:=\E\big[\ip{D(\Xi_t^V\psi(U_t))}{v}\big]=\E\Big[\Xi_t^V\psi(U_t)\int_0^tv(s)\cdot\dd W_s\Big],
\end{align*}
and
\begin{align*}
I_2:=\E\Big[\Xi_t^V\psi(U_t)\int_0^t\nabla V(U_s)\rho(s)\dd s\Big],\qquad
I_3:=\E\big[\Xi_t^V\nabla\psi(U_t)\rho(t)\big].
\end{align*}

\emph{Step 2: Estimate of the control term.}
Setting $K:=\lfloor t/2\rfloor$, we obtain from Lemma~\ref{lem:residual}$(iii)$ that
\[
\int_0^tv(s)\cdot\dd W_s=\sum_{0\le k\le(t-1)/2}\int_{2k}^{2k+1}v(s)\cdot\dd W_s+\mathbf 1_{t\in(2K,2K+1)}\int_{2K}^tv(s)\cdot\dd W_s .
\]
On a complete block, the Skorokhod integral is $\mathcal F_{2k+1}$-measurable, which allows us to condition at $2k+1$.
Proposition~\ref{prop:growth1}$(i)$ and monotonicity bound $|\PV{t-2k-1}\psi|$ by $C\norm\psi_\infty\fw_m\norm{\PV t\mathbf1}_{R_0}$.
Writing
\[
I_{1,k}:=\E\Big[\Xi_t^V\psi(U_t)\int_{2k}^{2k+1}v\cdot\dd W\Big]=\E\Big[\Xi^V_{2k+1}\Big(\int_{2k}^{2k+1}v\cdot\dd W\Big)(\PV{t-2k-1}\psi)(U_{2k+1})\Big],
\]
we apply Cauchy-Schwarz and Lemma~\ref{lem:residual}$(ii)$ with $n=2k$.
By Lemma~\ref{lem:poly}$(iii)$ and the bound $1+x^{4m}\le C_{m,\eta}e^{\eta x^2}$, we have $\E\fw_m(U_{2k+1})^2\le C_{m,\eta}e^{\eta\norm U^2}$. Hence
\begin{align*}
|I_{1,k}|&\le C\norm\psi_\infty\norm{\PV t\mathbf 1}_{R_0}e^{(2k+1)\bar\gamma}\big(\E\fw_m(U_{2k+1})^2\big)^{1/2}\Big(\E\Big|\int_{2k}^{2k+1}v\cdot\dd W\Big|^2\Big)^{1/2}\\
&\le C\norm\psi_\infty\norm{\PV t\mathbf 1}_{R_0}e^{\eta\norm U^2}e^{(2k+1)\bar\gamma-k\mathfrak a},
\end{align*}
and the series is summable because $\mathfrak a\ge2\bar\gamma+1$.
For the incomplete final block, we have $t<2K+1$, so Lemma~\ref{lem:residual}$(ii)$ gives
\begin{align*}
\Big|\E\Big[\Xi_t^V\psi(U_t)\int_{2K}^tv\cdot\dd W\Big]\Big|
&\le\norm\psi_\infty e^{t\bar\gamma}\Big(\E\Big|\int_{2K}^tv\cdot\dd W\Big|^2\Big)^{1/2}\\
&\le C\norm\psi_\infty e^{(2K+1)\bar\gamma-K\mathfrak a}e^{\eta\norm U^2}\le C\norm\psi_\infty e^{\eta\norm U^2}.
\end{align*}
Since $\norm{\PV t\mathbf 1}_{R_0}\ge1$, we obtain $|I_1|\le C\norm\psi_\infty\norm{\PV t\mathbf 1}_{R_0}e^{\eta\norm U^2}$.

\emph{Step 3: Estimates of the residual terms.}
For $I_2$ and $0\le s\le\lfloor t\rfloor$, Lemma~\ref{lem:FGRT31} gives $\rho(s)\in\mathcal F_{\lceil s\rceil}$. Conditioning at $\lceil s\rceil$ and applying Proposition~\ref{prop:growth1}$(i)$, we get
\begin{align*}
&~~\big|\E[\Xi_t^V\psi(U_t)\nabla V(U_s)\rho(s)]\big|\\
&\le C\norm\psi_\infty\norm{\nabla V}_\infty\norm{\PV t\mathbf1}_{R_0}
 e^{\lceil s\rceil\bar\gamma}\E[\norm{\rho(s)}\fw_m(U_{\lceil s\rceil})]\\
&\le C\norm\psi_\infty\norm{\nabla V}_\infty\norm{\PV t\mathbf1}_{R_0}
 e^{\eta\norm U^2+\lceil s\rceil\bar\gamma-\mathfrak a s/4}.
\end{align*}
The last bound follows from H\"older's inequality, Lemma~\ref{lem:residual}$(i)$ and Lemma~\ref{lem:poly}$(iii)$, with the exponents of the initial-state weights chosen to sum to $\eta$. Since $\lceil s\rceil\bar\gamma-\mathfrak a s/4\le\bar\gamma-(\mathfrak a/4-\bar\gamma)s$ and $\mathfrak a/4>\bar\gamma$, the time factor is integrable.
On the remaining interval $(\lfloor t\rfloor,t]$, whose length is at most one, we use $\Xi_t^V\le e^{t\bar\gamma}$ and Lemma~\ref{lem:residual}$(i)$.
This contribution is bounded by $C\norm\psi_\infty\norm{\nabla V}_\infty e^{\eta\norm U^2}$, and hence
\[
|I_2|\le C\norm\psi_\infty\norm{\nabla V}_\infty\norm{\PV t\mathbf1}_{R_0}e^{\eta\norm U^2}.
\]

Finally, $e^{\bar\gamma-\mathfrak a/4}\le\varepsilon_*$ gives
\[
|I_3|\le\norm{\nabla\psi}_\infty e^{t\bar\gamma}\E\norm{\rho(t)}\le C\norm{\nabla\psi}_\infty e^{\eta\norm U^2}\big(e^{\bar\gamma-\mathfrak a/4}\big)^t\le C\norm{\nabla\psi}_\infty e^{\eta\norm U^2}\varepsilon_*^t .
\]
Consequently,
\[
|I_1|+|I_2|+|I_3|\le
Ce^{\eta\norm U^2}\norm{\PV t\mathbf1}_{R_0}
\Big[(1+\norm{\nabla V}_\infty)\norm\psi_\infty
+\varepsilon_*^t\norm{\nabla\psi}_\infty\Big],
\]
which proves the assertion.
\end{proof}
The gradient estimate yields the uniform Feller property on each compact set $X_R$.

\begin{corollary}\label{cor:feller}
Assume (A1)-(A3), let $V,\psi\in C_b^1(H)$ and take $R\ge R_0(V)$.
The family $\{\norm{\PV t\mathbf1}_R^{-1}\PV t\psi\}_{t\ge0}$ is uniformly Lipschitz on $X_R$ in the metric of $H$.
In particular, it is uniformly equicontinuous.
\end{corollary}

\begin{proof}
As $X_R\subset B_H(R)$ is convex, the mean value inequality and Proposition~\ref{prop:grad} imply that $\PV t\psi$ is Lipschitz on $X_R$ with a constant not exceeding $C\norm{\PV t\mathbf1}_{R_0}$, where $C>0$ does not depend on $t$ (since $\varepsilon_*^t\le1$).
Since $X_{R_0}\subset X_R$, we have $\norm{\PV t\mathbf1}_{R_0}\le\norm{\PV t\mathbf1}_R$, which gives the required bound.
\end{proof}
\section{Proofs of the main results}\label{sec:proofs}
This section begins with the proof of multiplicative ergodicity. We then develop the pressure and equilibrium state properties needed for the large deviation principle. Finally, we prove Theorem~\ref{thm:A} and Corollary~\ref{cor:level1}, including the uniqueness of the zero and the non-triviality of the rate function.

\subsection{\texorpdfstring{Proof of Theorem~\ref{thm:B}}{Proof of Theorem 2.3}}\label{subsec:proofB}
We begin by approximating bounded continuous functions by smooth functions on compact sets. This enables us to extend the semigroup convergence and the pressure limits.

\begin{lemma}\label{lem:buc}
Let $f\in C_b(H)$.
There exist $f_n\in C_b^1(H)$ with $\norm{f_n}_\infty\le\norm f_\infty$.
For every compact $K\subset H$, they satisfy $\sup_K|f_n-f|\to0$ as $n\to\infty$.
\end{lemma}

\begin{proof}
Let $\Pi_n$ be the weighted orthogonal projection onto $H_n:=\spn\{\sigma_k^l,\psi_k^l:|k|\le n\}$.
Identifying $H_n$ with $\R^{d_n}$ through an orthonormal basis and taking a standard mollifier $\varphi_\varepsilon$, we set
\[
f_n(u):=\int_{H_n}f(\Pi_nu-z)\varphi_{1/n}(z)\dd z.
\]
Then $f_n\in C_b^1(H)$ and $\norm{f_n}_\infty\le\norm f_\infty$.
Since $\Pi_n$ converges to the identity uniformly on any compact set $K\subset H$, the points $\Pi_nu-z$ with $n\ge1$, $u\in K$ and $\norm z\le1/n$ form a relatively compact set, and $\norm{\Pi_nu-z-u}\to0$ as $n\to\infty$, uniformly in such $u$ and $z$.
Hence the uniform continuity of $f$ on the closure of this set gives $\sup_K|f_n-f|\to0$.
\end{proof}

\begin{proof}[Proof of Theorem~\ref{thm:B}]
Let $\tilde V:=V-\inf V\ge0$.
Since $\Osc(\tilde V)=\Osc(V)$, the constants $(m,R_0)$ in Proposition~\ref{prop:growth1} are the same for $V$ and $\tilde V$.

\emph{Step 1: Application of the abstract theorem.}
We verify its hypotheses on $H$ with compact sets $X_R$ and weight $\fw_m$, using the results of Sections~\ref{sec:apriori} and~\ref{sec:fk}.
\begin{enumerate}
\item[(a)] Corollary~\ref{cor:expmom}$(ii)$ gives continuity of $U\mapsto P_t^{\tilde V}(U,\cdot)$ in the weak topology.
The kernels are non-zero because $P_t^{\tilde V}\mathbf1\ge1$, and they satisfy the Chapman-Kolmogorov relation by the Markov property.
For $t>0$, they are concentrated on $X_\infty=H^1$ by Proposition~\ref{prop:wp}$(i)$ and absolute continuity with respect to $P_t$.
\item[(b)] The compact sets $X_R$ are increasing and have dense union by Subsection~\ref{subsec:irr}.
Proposition~\ref{prop:growth1}$(ii)$, $(iii)$ gives the growth conditions
\[
\sup_{t\ge0}\frac{\norm{P_t^{\tilde V}\fw_m}_{L^\infty_{\fw_m}}}
{\norm{P_t^{\tilde V}\mathbf1}_{R_0}}<\infty,
\qquad
\sup_{t\in[0,1]}\norm{P_t^{\tilde V}\mathbf1}_\infty<\infty.
\]
\item[(c)] Lemma~\ref{lem:timecont} gives time continuity of $P_t^{\tilde V}g(U)$ for $g\in C_{\fw_m}(H)$ and $U\in H$.
\item[(d)] Proposition~\ref{prop:irr} gives uniform irreducibility for every $r_*,R\ge1$ and $r>0$.
Here $l=l(r_*,r,R)$ may depend on the initial radius, while $p=p(r_*,r)>0$ is independent of it, as required.
\item[(e)] For $\psi\in C_b^1(H)$ and $R\ge R_0$, Corollary~\ref{cor:feller} gives uniform equicontinuity of $\norm{P_t^{\tilde V}\mathbf1}_R^{-1}P_t^{\tilde V}\psi$ on $X_R$.
\item[(f)] Proposition~\ref{prop:eigen} supplies an eigenmeasure $\mu_{t,\tilde V}\in\cP_{\fw_m}(H)$ for each $t>0$.
It satisfies the tail condition
\[
\norm{P_t^{\tilde V}\fw_m}_R\int_{H\setminus X_R}\fw_m\dd\mu_{t,\tilde V}\longrightarrow0
\qquad\text{as }R\to\infty.
\]
Its eigenrelation and (a) imply $\mu_{t,\tilde V}(X_\infty)=1$.
\end{enumerate}
The measure-determining class $\mathcal C=C_b^1(H)$ contains $\mathbf1$.
Indeed, it contains $\cos\ip\cdot\xi$ and $\sin\ip\cdot\xi$ for every $\xi\in H$, so equality on $\mathcal C$ implies equality of characteristic functionals and hence of measures.

Properties (a)-(f) verify the assumptions of~\cite[Theorem~6.6]{Ner19}, equivalently~\cite[Theorem~7.4]{MN18a}.
For time continuity, the weighted closure $\mathcal C_w$ is contained in $C_{\fw_m}(H)$ because $\fw_m$ is bounded on bounded sets, so (c) applies.
The abstract theorem therefore gives the following conclusions.
\begin{enumerate}
\item[(C1)] For each $t>0$, $\mu_{t,\tilde V}$ is the unique eigenmeasure in $\cP_{\fw_m}(H)$ satisfying (f) and concentrated on $X_\infty$.
The measure $\mu:=\mu_{t,\tilde V}$ is independent of $t$.
Moreover, $\lambda_{t,\tilde V}=\lambda_{\tilde V}^t$ for some $\lambda_{\tilde V}>0$, and $\supp\mu=H$.
\item[(C2)] There is a non-negative $h\in L^\infty_{\fw_m}(H)$ with $\ip h\mu=1$ and $P_t^{\tilde V}h=\lambda_{\tilde V}^th$ for $t>0$.
Its restriction to each $X_R$ is continuous.
For $\psi\in C_b^1(H)$ and $R\ge1$ we have
\begin{equation}\label{eq:abs-conv}
\lambda_{\tilde V}^{-t}P_t^{\tilde V}\psi\to\ip\psi\mu h\quad\text{in }C(X_R)\cap L^1(H,\mu)\quad(t\to\infty).
\end{equation}
\item[(C3)] If a Borel set $B\subset H$ satisfies the following condition for some $s>0$
\[
\sup_{U\in B}\int_{H\setminus X_R}\fw_m(v)\,P_s^{\tilde V}(U,\dd v)\to0\quad(R\to\infty),
\]
then $\lambda_{\tilde V}^{-t}P_t^{\tilde V}\psi\to\ip\psi\mu\,h$ in $L^\infty(B)$ for $\psi\in C_b^1(H)$.
\end{enumerate}

By Proposition~\ref{prop:eigen}$(ii)$, $\mu\in\bigcap_{\varkappa\in(0,\eta_1)}\cP_{\fm_\varkappa}(H)$.
Returning to the original potential, we set $\lambda:=e^{\inf V}\lambda_{\tilde V}$, so that $\lambda_{\tilde V}^{-t}P^{\tilde V}_t=\lambda^{-t}\PV t$ and $\PV th=\lambda^th$ on $H$.
We henceforth write $\lambda_V:=\lambda$, $h_V:=h$ and $\mu_V:=\mu$.

\emph{Step 2: Convergence for bounded continuous test functions.}
We first extend the convergence for $\psi\in C_b^1(H)$ from $X_R$ to $B_H(R)$.
By (C3), it suffices to verify the tail condition with $s=1$.
Lemma~\ref{lem:H1mom} with $T=2$, $p=2m+2$ and $\eta=1$ gives precisely the required estimate:
\[
\sup_{\norm U\le R}\int_{H\setminus X_{r_*}}\fw_m(v)P^{\tilde V}_1(U,\dd v)\le C(m)e^{\Osc(V)+R^2}r_*^{-2}\to0\quad(r_*\to\infty).
\]
Thus $\lambda^{-t}\PV t\psi$ converges uniformly on $B_H(R)$ for $\psi\in C_b^1(H)$.
Taking $\psi=\mathbf1$, we obtain continuity of $h$ on $H$, and hence $h\in C_{\fw_m}(H)$ since $h\in L^\infty_{\fw_m}(H)$.
The same convergence, together with boundedness on finite time intervals, gives
\begin{equation}\label{eq:lambda-bound}
\sup_{t\ge0}\lambda^{-t}\norm{\PV t\mathbf 1}_{R_0}<\infty .
\end{equation}
To pass to $\psi\in C_b(H)$, we combine the compact approximation in Lemma~\ref{lem:buc} with tail bounds.
For $0<\varkappa<\varkappa'<\eta_1$, Propositions~\ref{prop:growth1} and~\ref{prop:growth23}, together with \eqref{eq:lambda-bound}, give
\begin{equation}\label{eq:tail-bounds}
\begin{gathered}
\lambda^{-t}\PV t\mathbf 1(U)\le C\fw_m(U),\qquad \lambda^{-t}\PV t\fm_{\varkappa'}(U)\le C\fm_{\varkappa'}(U)\qquad(t\ge0),\\
\lambda^{-t}\PV t\mathsf F(U)\le C\fm_\varkappa(U)\qquad(t\ge1).
\end{gathered}
\end{equation}
Here $\mathsf F$ is the extended $H^1$ energy from Proposition~\ref{prop:growth23}.
Proposition~\ref{prop:eigen}$(ii)$ also gives $\ip{\mathsf F}\mu+\ip{\fm_{\varkappa'}}\mu<\infty$.

For $\psi\in C_b(H)$, choose $\psi_n$ from Lemma~\ref{lem:buc}.
Set $\delta_n(r):=\sup_{X_r}|\psi-\psi_n|$, which tends to zero for fixed $r$.
Inserting $|\psi-\psi_n|\le\delta_n(r)+2\norm\psi_\infty r^{-2}\mathsf F$ into \eqref{eq:tail-bounds}, for $U\in B_H(R)$ and $t\ge1$, we obtain
\[
\lambda^{-t}\big|\PV t(\psi-\psi_n)(U)\big|+\big|\ip{\psi-\psi_n}\mu\big|\,h(U)\le C_R\big(\delta_n(r)+r^{-2}\big),
\]
where $C_R$ depends on $R,\norm\psi_\infty$, but not on $t,n,r$.
Choosing first $r$ and then $n$ sufficiently large makes this error arbitrarily small, uniformly in $t$.
Combined with convergence for $\psi_n$, this proves uniform convergence for $\psi$ on $B_H(R)$.

\emph{Step 3: Convergence for exponentially growing test functions.}
For $\psi\in C_{\fm_\varkappa}(H)$, we use the truncations $\psi_N:=(\psi\wedge N)\vee(-N)$.
The zero function requires no argument, so assume $\psi\ne0$ and set $A_N:=\varkappa^{-1}\log(N/\norm\psi_{L^\infty_{\fm_\varkappa}})$.
Then
\[
|\psi-\psi_N|\le\norm\psi_{L^\infty_{\fm_\varkappa}}\fm_\varkappa\mathbf 1_{\{\norm U^2>A_N\}}\le\norm\psi_{L^\infty_{\fm_\varkappa}}e^{-(\varkappa'-\varkappa)A_N}\fm_{\varkappa'} .
\]
By \eqref{eq:tail-bounds}, the semigroup and eigenmeasure errors on $B_H(R)$ are bounded by $C_Re^{-(\varkappa'-\varkappa)A_N}$, which tends to zero uniformly in $t$.
The approximation argument of Step 2 therefore gives uniform convergence for $\psi\in C_{\fm_\varkappa}(H)$.

The same approximations give convergence in $L^1(H,\mu)$.
Indeed, the eigenrelation yields $\int_H|\lambda^{-t}\PV t\phi|\dd\mu\le\ip{|\phi|}\mu$, and $\norm h_{L^1(\mu)}=1$.
Dual convergence follows by dominated convergence using the first bound in \eqref{eq:tail-bounds}.

\emph{Step 4: Uniformity over the initial distributions.}
Let $\varkappa'>\varkappa$ and $M>0$ be the numbers in $(iii)$.
By Proposition~\ref{prop:growth23}$(i)$ and \eqref{eq:lambda-bound}, we have $|\lambda^{-t}\PV t\psi|\le C\fm_\varkappa$ for $\psi\in C_{\fm_\varkappa}(H)$.
On the other hand,
\[
\sup_{\nu\in\Lam(\varkappa',M)}
\int_{H\setminus B_H(R)}\fm_\varkappa\dd\nu
\le Me^{-(\varkappa'-\varkappa)R^2}.
\]
Combining these two bounds with the convergence on $B_H(R)$, we obtain $(iii)$.

\emph{Step 5: Positivity and uniqueness.}
Since $h$ is continuous, $\ip h\mu=1$ and $H^1$ is dense in $H$, there are $\hat U\in H^1$ and $r,c>0$ such that $h>c$ on $B_H(\hat U,r)$.
Let $r_*\ge\max\{1,\norm{\hat U}_{H^1}\}$.
By the eigenrelation and Proposition~\ref{prop:irr}, for $U\in X_R$ with $R\ge1$, we have $h(U)\ge\lambda_{\tilde V}^{-l}c\,p(r_*,r)>0$, where $l=l(r_*,r,R)$.
Hence $h>0$ on $H^1$.
As $U_1\in H^1$ almost surely, the eigenrelation gives $h(U)\ge\lambda^{-1}e^{-\norm V_\infty}\E_Uh(U_1)>0$ for any $U\in H$.

Now we turn to the uniqueness.
Any eigenmeasure $\mu'\in\cP_{\fm_\varkappa}(H)$ satisfies the moment and tail bounds of Proposition~\ref{prop:eigen} and, by its eigenrelation and (a), is concentrated on $X_\infty$.
Hence $\mu'=\mu$ by (C1).
If $h'\in C_{\fw_m}(H)$ satisfies $\PV th'=\lambda'^th'$ and $\ip{h'}\mu=1$, then integrating this relation against $\mu$, we get $\lambda'=\lambda$.
Since $C_{\fw_m}(H)\subset C_{\fm_\varkappa}(H)$, the convergence in $(ii)$ gives
$h'=\lambda^{-t}\PV th'\to\ip{h'}\mu h=h$ as $t\to\infty$.
\end{proof}

\subsection{Exponential tightness, pressure and equilibrium states}\label{sec:ldp}\label{subsec:equil}
The Kifer-type criterion in~\cite[Theorem~3.3]{JNPS18} requires exponential tightness, pressure limits and unique equilibrium states for a dense family of potentials.
We establish these properties in that order.
The equilibrium-state argument uses a Doob transform.
The integrated energy estimate provides a compactness function for exponential tightness.
\begin{lemma}\label{lem:tight}
Assume (A1)-(A3), and let $\gamma\in(0,\eta_1]$ and $M\ge1$.
Define $\mathsf E(U):=\frac{\gamma\kappa\norm U_{H^1}^2}{2}$ on $H^1$ and $\mathsf E(U):=+\infty$ otherwise.
Then,
\begin{enumerate}
\item[(i)] $\mathsf E$ is lower semicontinuous on $H$, and for every $\alpha\ge0$ the level set $\{\mathsf E\le\alpha\}=X_{(2\alpha/(\gamma\kappa))^{1/2}}$ is compact in $H$;
\item[(ii)] for any $\nu\in\Lam(\gamma,M)$ and $t>0$,
\[
\E_\nu\exp\big(t\ip{\mathsf E}{L_t}\big)=\E_\nu\exp\Big(\int_0^t\mathsf E(U_s)\dd s\Big)\le Me^{\gamma\fB_0t}.
\]
\end{enumerate}
\end{lemma}

\begin{proof}
For $(i)$, the lower semicontinuity of the extended $H^1$ norm and the compact embedding $H^1\hookrightarrow\hookrightarrow H$ give the assertion.
Integrating Lemma~\ref{lem:expmom}$(ii)$, with $\eta=\gamma$, against $\nu$ yields $(ii)$.
In particular, for $A>0$,
\[
\mathcal K_A:=\big\{\sigma\in\cP(H):\ip{\mathsf E}{\sigma}\le A\big\}
\quad\text{is compact in }\cP(H),
\]
and Markov's inequality gives
\begin{equation}\label{eq:exp-tight-prob}
\sup_{\nu\in\Lam(\gamma,M)}\Prob_\nu\{L_t\notin\mathcal K_A\}
\le e^{-At}\sup_{\nu\in\Lam(\gamma,M)}
\E_\nu e^{t\ip{\mathsf E}{L_t}}
\le M e^{-t(A-\gamma\fB_0)}.
\end{equation}
\end{proof}
Multiplicative ergodicity gives pressure limits for smooth potentials.
The preceding tightness estimate allows us to extend them to all bounded continuous potentials.
\begin{lemma}\label{lem:pressure}
Assume (A1)-(A3), and let $\gamma\in(0,\eta_1]$, $M\ge1$ and $\cV:=C_b^1(H)$. For $V\in C_b(H)$, $\nu\in\cP(H)$ and $t>0$, set
\[
Q_t^\nu(V):=\frac1t\log\E_\nu\exp\Big(\int_0^tV(U_s)\dd s\Big).
\]
Then:
\begin{enumerate}
\item[(i)] for any $V\in\cV$, $Q_t^\nu(V)\to\log\lambda_V$ as $t\to\infty$, uniformly with respect to $\nu\in\Lam(\gamma,M)$;
\item[(ii)] for every $V\in C_b(H)$, $Q_t^\nu(V)$ converges to a limit $Q(V)$ uniformly over $\nu\in\Lam(\gamma,M)$.
The limit is independent of $\nu$;
\item[(iii)] $Q$ is a $1$-Lipschitz convex function on $C_b(H)$, $Q(0)=0$, and $Q(V+c)=Q(V)+c$ for $V\in C_b(H)$ and $c\in\R$.
\end{enumerate}
\end{lemma}

\begin{proof}
For $(i)$, fix $\varkappa\in(0,\gamma)$ and apply Theorem~\ref{thm:B}$(iii)$ with $\psi=\mathbf1$ and $\varkappa'=\gamma$.
It then suffices to bound $\ip{h_V}\nu$ above and away from zero uniformly over $\Lam(\gamma,M)$.
The upper bound follows from $h_V\le C\fw_m$.
For the lower bound, choose $R$ so that $e^{\gamma R^2}\ge2M$, which ensures $\nu(B_H(R))\ge1/2$.
Lemma~\ref{lem:H1mom} gives $r_*=r_*(R)$ such that $\Prob_U(U_1\in X_{r_*})\ge1/2$ for all $U\in B_H(R)$.
The eigenrelation now implies
\[
h_V(U)=\lambda_V^{-1}\PV 1h_V(U)\ge\lambda_V^{-1}e^{-\norm V_\infty}\Big(\inf_{X_{r_*}}h_V\Big)\Prob_U(U_1\in X_{r_*})\ge\frac{e^{-\norm V_\infty}}{2\lambda_V}\inf_{X_{r_*}}h_V=:c_R,
\]
where $c_R>0$ by compactness and positivity of $h_V$.
Thus $\ip{h_V}\nu\ge c_R/2$, and
\[
Q_t^\nu(V)-\log\lambda_V=\frac1t\log\big(\lambda_V^{-t}\E_\nu\Xi_t^V\big)\to0\qquad(t\to\infty)
\]
uniformly with respect to $\nu\in\Lam(\gamma,M)$.

For $(ii)$, we use the compact approximation argument of~\cite[Section~5.6]{JNPS18}, based on~\cite[Proposition~3.17]{FK06}.
To make the argument uniform over the initial measures, choose $V_n\in\cV$ from Lemma~\ref{lem:buc} and set $a:=\sup_n\norm{V_n}_\infty\vee\norm V_\infty$.
On the compact set $K_\alpha:=\{\mathsf E\le\alpha\}$, the error $\delta_n(\alpha):=\sup_{K_\alpha}|V-V_n|$ tends to zero for each fixed $\alpha$.
Choose $p\in(1,2]$, set $q:=p/(p-1)$ and take $\alpha\ge2aq$.
For $W:=\pm q(V-V_n)$, the bound $\mathbf1_{H\setminus K_\alpha}\le\mathsf E/\alpha$ gives
$\ip{|W|}{L_t}\le q\delta_n(\alpha)+(2aq/\alpha)\ip{\mathsf E}{L_t}$.
Lemma~\ref{lem:tight}$(ii)$ and H\"older's inequality therefore yield
\[
\frac1qQ_t^\nu(W)\le\delta_n(\alpha)+\frac{2a}\alpha\Big(\gamma\fB_0+\frac{\log M}t\Big)=:e_n(\alpha,t).
\]
Each $Q_t^\nu$ is convex by H\"older's inequality and is $1$-Lipschitz in the supremum norm.
Applying these properties to $V=V_n+(V-V_n)$ and $V_n=V+(V_n-V)$, with conjugate exponents $p,q$, we obtain
\[
pQ_t^\nu(V_n)-(p-1)a-p\,e_n(\alpha,t)\le Q_t^\nu(V)\le\frac1p\big(Q_t^\nu(V_n)+(p-1)a\big)+e_n(\alpha,t).
\]
By $(i)$, $Q_t^\nu(V_n)\to\log\lambda_{V_n}$ uniformly over $\Lam(\gamma,M)$, and $|\log\lambda_{V_n}|\le a$.
It follows that
\begin{align*}
&\limsup_{t\to\infty}\sup_{\nu\in\Lam(\gamma,M)}Q_t^\nu(V)
-\liminf_{t\to\infty}\inf_{\nu\in\Lam(\gamma,M)}Q_t^\nu(V)\\
&\quad\le 2a\Big(p-\frac1p\Big)
+(p+1)\Big(\delta_n(\alpha)+\frac{2a\gamma\fB_0}{\alpha}\Big).
\end{align*}
Given $\varepsilon>0$, choose $p\in(1,2]$ so close to $1$ that
$2a(p-p^{-1})<\varepsilon/3$, and set $q=p/(p-1)$.
Next choose $\alpha\ge\max\{1,2aq\}$ such that
$2a\gamma\fB_0(p+1)/\alpha<\varepsilon/3$.
For this fixed $\alpha$, choose $n$ so large that
$(p+1)\delta_n(\alpha)<\varepsilon/3$.
The preceding difference of the uniform upper and lower limits is then at most $\varepsilon$.
Since $\varepsilon$ is arbitrary, $(ii)$ follows.

The remaining properties follow directly from $Q_t^\nu(0)=0$ and $Q_t^\nu(V+c)=Q_t^\nu(V)+c$.
Convexity and the Lipschitz bound also pass to the limit, proving $(iii)$.
\end{proof}

For $\sigma\in\cP(H)$, let
\[
I(\sigma):=\sup_{V\in C_b(H)}\big(\ip V\sigma-Q(V)\big)
\]
be the Legendre transform of the pressure function $Q$ of Lemma~\ref{lem:pressure}$(ii)$.
We next identify the unique equilibrium state for each $V\in\cV$.
\begin{lemma}\label{lem:equil}
Assume (A1)-(A3), and let $V\in\cV$ and $\sigma_V:=h_V\mu_V\in\cP(H)$. Then, 
\begin{enumerate}
\item[(i)] $I(\sigma_V)=\ip V{\sigma_V}-Q(V)$;
\item[(ii)] if $\sigma\in\cP(H)$ satisfies $I(\sigma)=\ip V\sigma-Q(V)$, then $\sigma=\sigma_V$.
\end{enumerate}
\end{lemma}

\begin{proof}
Since $I$ is the Legendre transform of $Q$, a measure $\sigma\in\cP(H)$ satisfies $I(\sigma)=\ip V\sigma-Q(V)$ if and only if
\begin{equation}\label{eq:equil-ineq}
\ip f\sigma\le Q(V+f)-Q(V),\qquad f\in C_b(H).
\end{equation}
The exponential moment bound in Lemma~\ref{lem:tight}$(ii)$ and the compact sublevel sets of $\mathsf E$ verify~\cite[Lemma~3.2]{JNPS18}.
Using the pressure limit of Lemma~\ref{lem:pressure}, Step~1 of the proof of~\cite[Theorem~3.3]{JNPS18} gives goodness of $I$ and existence of a measure satisfying \eqref{eq:equil-ineq}.
It therefore suffices to prove $(ii)$.

We follow the Doob transform argument of~\cite[Proposition~1.3]{MN18a}.
Related formulations appear in~\cite[Section~5.7]{JNPS18} and~\cite[Section~2]{Ner19}.
For $f,g\in C_b(H)$, set
\[
\mathfrak S_t^fg:=\lambda_V^{-t}h_V^{-1}P_t^{V+f}(h_Vg),
\qquad S_t:=\mathfrak S_t^0.
\]
Corollary~\ref{cor:expmom}$(ii)$, Lemma~\ref{lem:timecont} and the eigenrelation show that $(S_t)$ is a Markov semigroup on $C_b(H)$, pointwise continuous in time.
Theorem~\ref{thm:B}$(ii)$, applied to $h_Vg$, gives $S_tg\to\ip g{\sigma_V}$ pointwise.
For $\psi\in C_b(H)$ with $\inf\psi>0$ and $s>0$, set
\[
\bar\psi:=\frac1s\int_0^sS_\tau\psi\dd\tau,\qquad \mathbb L\bar\psi:=\frac{S_s\psi-\psi}s,\qquad f:=-\frac{\mathbb L\bar\psi}{\bar\psi},\qquad V_1:=V+f.
\]
Then $\bar\psi,f\in C_b(H)$ and $\inf\bar\psi\ge\inf\psi$.
The semigroup identity gives $S_t\bar\psi-\bar\psi=\int_0^tS_r(\mathbb L\bar\psi)\dd r$.
Since $f\bar\psi=-\mathbb L\bar\psi$, the bounded-potential Duhamel formula~\cite[Proposition~1.3]{MN18a} gives
\[
\mathfrak S_t^f\bar\psi-\bar\psi
=\int_0^t\mathfrak S_{t-r}^f(\mathbb L\bar\psi+f\bar\psi)\dd r=0.
\]
Consequently, we obtain
\[
P_t^{V_1}(h_V\bar\psi)=\lambda_V^th_V\bar\psi,\qquad t\ge0 .
\]
The remaining point is to prove $Q(V_1)\le Q(V)$ with the weights available here.
Fix $\varkappa\in(0,\eta_1)$ and take $(m',R_0')$ from the growth estimates for $V_1$.
Set $c_r:=\inf_{X_r}(h_V\bar\psi)>0$.
Splitting the terminal state according to membership in $X_r$, we use the eigenrelation on $X_r$ and Proposition~\ref{prop:growth23}$(ii)$ on its complement.
For $U\in X_{R_0'}$ and $t\ge1$, the identity
\[
\mathbf1\le c_r^{-1}h_V\bar\psi+\mathbf1_{H\setminus X_r}
\le c_r^{-1}h_V\bar\psi+r^{-2}\mathsf F
\]
and Proposition~\ref{prop:growth23}$(ii)$ give
\begin{align*}
P_t^{V_1}\mathbf1(U)
&\le c_r^{-1}\lambda_V^t h_V(U)\bar\psi(U)
+r^{-2}P_t^{V_1}\mathsf F(U)\\
&\le c_r^{-1}\lambda_V^t\sup_{X_{R_0'}}(h_V\bar\psi)
+Cr^{-2}e^{\varkappa R_0'^2}\norm{P_t^{V_1}\mathbf1}_{R_0'}.
\end{align*}
The constant $C$ is independent of $r$ and $t\ge1$.
Fix $r$ sufficiently large that $Cr^{-2}e^{\varkappa R_0'^2}\le1/2$.
Taking the supremum over $X_{R_0'}$ yields
\[
\norm{P_t^{V_1}\mathbf1}_{R_0'}
\le2c_r^{-1}\lambda_V^t\sup_{X_{R_0'}}(h_V\bar\psi).
\]
Since $0\in X_{R_0'}$, Lemma~\ref{lem:pressure} gives
\[
Q(V_1)=\lim_{t\to\infty}t^{-1}\log P_t^{V_1}\mathbf1(0)
\le\log\lambda_V=Q(V).
\]
Equation \eqref{eq:equil-ineq} now implies $\ip{\mathbb L\bar\psi/\bar\psi}\sigma\ge0$.

As in~\cite[Proposition~1.3]{MN18a}, taking $\psi=1+\vartheta g$ with $g\in C_b(H)$ and sufficiently small $|\vartheta|$, dividing the last inequality by $\vartheta$ and letting $\vartheta\to0$ from both sides, we derive from the dominated convergence theorem that $\ip{S_sg-g}\sigma=0$ for any $s>0$.
Hence $\sigma$ is invariant under $S_t$.
Since $S_tg\to\ip g{\sigma_V}$ and $|S_tg|\le\norm g_\infty$, applying the dominated convergence theorem again, we get $\ip g\sigma=\ip g{\sigma_V}$ for any $g\in C_b(H)$, and hence $\sigma=\sigma_V$.
\end{proof}

\subsection{\texorpdfstring{Proofs of Theorem~\ref{thm:A} and Corollary~\ref{cor:level1}}{Proofs of Theorem 2.4 and Corollary 2.5}}\label{subsec:proofA}
We first combine the results of Subsection~\ref{subsec:equil} with the Kifer-type criterion to establish parts $(i)$ and $(iv)$ of Theorem~\ref{thm:A}, which will then be used to analyze the zero and non-triviality of the rate function.

\begin{proof}[Proof of Theorem~\ref{thm:A}(i) and (iv)]
Let $M\ge1$.
Since $\Lam(\gamma,M)\subset\Lam(\gamma',M)$ for $\gamma'\le\gamma$, it suffices to consider the case $0<\gamma\le\eta_1$.
We consider the directed set
\[
\mathcal I:=(0,\infty)\times\Lam(\gamma,M),\qquad
(t,\nu)\preceq(t',\nu')\ \Longleftrightarrow\ t<t'\ \text{or}\ (t,\nu)=(t',\nu'),
\]
and, for $\iota=(t,\nu)$, we set $r_\iota:=t$ and denote by $\mathsf L_\iota$ the random measure $L_t$ under $\Prob_\nu$.
As in~\cite{Ner19}, the limits along $\mathcal I$ are the limits as $t\to\infty$ which are uniform in $\nu\in\Lam(\gamma,M)$.

By Lemma~\ref{lem:pressure}, the pressure limit exists for any $V\in C_b(H)$.
The compactness of $\mathcal K_A$ and \eqref{eq:exp-tight-prob} give the uniform exponential tightness in $\cP(H)$.
Since the algebra $\cV=C_b^1(H)$ contains the constants and separates the points of $H$, its restriction to any compact set is dense by the Stone-Weierstrass theorem.
Moreover, by Lemma~\ref{lem:equil}, $h_V\mu_V$ is the unique equilibrium state for any $V\in\cV$.
Hence the hypotheses of~\cite[Theorem~3.3]{JNPS18} are satisfied, and we obtain the uniform large deviation bounds with the good rate function $I$.

The properties of $Q$ follow from Lemma~\ref{lem:pressure}$(iii)$.
Since $\delta_0\in\Lam(\gamma,M)$ for any $\gamma>0$ and $M\ge1$, the function $Q$ does not depend on $\gamma$ and $M$, and hence neither does its Legendre transform $I$.
\end{proof}

We next turn to assertions $(ii)$ and $(iii)$ of Theorem~\ref{thm:A}.
Non-triviality of the rate function follows by keeping trajectories near the origin at an exponential probability cost.
First we establish stability under a uniformly small driving path.

\begin{lemma}\label{lem:iss}
Assume (A1)-(A3).
There is $c_0>0$, depending only on $\nu_1,\nu_2,g,\sigma_\theta$, with the following property.
Let $\varepsilon\in(0,1]$ and $t>0$.
If $\norm{U_0}\le\varepsilon/2$ and $\zeta\in C_0([0,t];\R^4)$ satisfies $|\zeta|_\infty\le c_0\varepsilon$, then
\[
\norm{\Phi_s(U_0,\zeta)}\le\varepsilon,\qquad s\in[0,t].
\]
\end{lemma}

\begin{proof}
Set $\xi=\sigma_\theta\zeta$ and $\bar U=\Phi_s(U_0,\zeta)-\xi$, which solves \eqref{eq:shifted}. By the same energy computation as in Lemma~\ref{lem:phi}$(ii)$, we obtain
\[
\frac12\frac{\dd}{\dd s}\norm{\bar U}^2+Z(\bar U)=-\ip{B(\bar U+\xi)}{\bar U}_{\fc}+\ip{G(\bar U+\xi)}{\bar U}_{\fc}-\ip{A\xi}{\bar U}_{\fc}.
\]
The transport cancellation gives $\ip{B(\bar U+\xi,\bar U)}{\bar U}_\fc=0$.
Since $\xi_w=0$, we also have $B(\bar U+\xi,\xi)=B(\bar U,\xi)$.
Using \eqref{eq:buoy}, choose constants $\beta_1,\beta_2,\beta_3$ depending only on $\nu_1,\nu_2,g,\sigma_\theta$ such that
\begin{align*}
|\ip{B(\bar U,\xi)}{\bar U}_{\fc}|=|\ip{(K*\bar w)\cdot\nabla\xi_\theta}{\bar\theta}|\le\beta_1|\zeta|_\infty\norm{\bar U}^2,
\end{align*}
and 
\begin{align*}
\ip{G(\bar U+\xi)}{\bar U}_{\fc}\le\frac12Z(\bar U)+\beta_2|\zeta|_\infty\norm{\bar U},\qquad|\ip{A\xi}{\bar U}_{\fc}|\le\beta_3|\zeta|_\infty\norm{\bar U}.
\end{align*}
Setting $\beta_4:=4(\beta_2+\beta_3)^2/\kappa$, we use $Z(\bar U)\ge\kappa\norm{\bar U}^2$ and Young's inequality to absorb the terms linear in $\norm{\bar U}$.
For $|\zeta|_\infty\le\kappa/(8\beta_1)$, we obtain
\[
\frac{\dd}{\dd s}\norm{\bar U}^2\le-\kappa\norm{\bar U}^2+2\beta_1|\zeta|_\infty\norm{\bar U}^2+2(\beta_2+\beta_3)|\zeta|_\infty\norm{\bar U}\le-\frac\kappa2\norm{\bar U}^2+\beta_4|\zeta|_\infty^2.
\]
Gronwall's inequality gives $\norm{\bar U(s)}^2\le\norm{U_0}^2+(2\beta_4/\kappa)|\zeta|_\infty^2$.
Choose
\[
c_0:=\min\Big\{\frac\kappa{8\beta_1},
\frac1{2\big((2\beta_4/\kappa)^{1/2}+\norm{\sigma_\theta}_{L(\R^4,H)}\big)}\Big\}.
\]
Since $\varepsilon\le1$ and $|\zeta|_\infty\le c_0\varepsilon$, we get
\[
\norm{\Phi_s(U_0,\zeta)}\le\norm{U_0}
+\Big((2\beta_4/\kappa)^{1/2}+\norm{\sigma_\theta}_{L(\R^4,H)}\Big)|\zeta|_\infty\le\varepsilon.
\qedhere
\]
\end{proof}
We combine the small-ball estimate with the large deviation upper bound to obtain finite-rate measures supported near the origin.

\begin{proposition}\label{prop:nontriv}
Assume (A1)-(A3).
Let $\gamma\in(0,\eta_1)$ and $M\ge\max\{1,\ip{\fm_\gamma}{\mu_*}\}$.
This moment is finite by Corollary~\ref{cor:expmom}$(iv)$.
Suppose Theorem~\ref{thm:A}$(i)$ holds for $\Lam(\gamma,M)$.
Then:
\begin{enumerate}
\item[(i)] $\mu_*\in\Lam(\gamma,M)$, $\supp\mu_*=H$ and $I(\mu_*)=0$;
\item[(ii)] for every $\varepsilon\in(0,1]$, there exists $\sigma_\varepsilon\in\cP(H)$ such that
\[
\sigma_\varepsilon(B_H(0,\varepsilon))=1,
\qquad I(\sigma_\varepsilon)<\infty,
\qquad \sigma_\varepsilon\ne\mu_*.
\]
\end{enumerate}
\end{proposition}

\begin{proof}
By the choice of $M$, we have $\mu_*\in\Lam(\gamma,M)$.
Integrating the eigenmeasure relation in Theorem~\ref{thm:B}$(i)$ with $V=0$, we get $\lambda_0=1$.
Hence $\mu_0$ is invariant, and $\mu_0=\mu_*$ by Proposition~\ref{prop:wp}$(iii)$.
This implies that $\supp\mu_*=H$, since $\supp\mu_0=H$ by Theorem~\ref{thm:B}$(i)$.
Moreover, by the invariance of $\mu_*$ and Jensen's inequality, for any $V\in C_b(H)$,
\[
\E_{\mu_*}\exp\Big(\int_0^tV(U_s)\dd s\Big)\ge\exp\big(t\ip V{\mu_*}\big).
\]
Passing to the pressure limit, we get $Q(V)\ge\ip V{\mu_*}$, and hence $I(\mu_*)\le0$.
Since $I\ge0$, this proves $(i)$.

For $(ii)$, fix $\varepsilon\in(0,1]$ and set
\[
F_\varepsilon:=\{\sigma\in\cP(H):\sigma(B_H(0,\varepsilon))=1\},
\]
which is closed by the portmanteau theorem.
Since $\supp\mu_*=H$, we have $\mu_*(H\setminus B_H(0,\varepsilon))>0$, and hence $\mu_*\notin F_\varepsilon$.
Set $E_t:=\{\sup_{s\in[0,t]}|W_s|\le c_0\varepsilon\}$.
Starting from $U_0=0$, Lemmas~\ref{lem:phi} and~\ref{lem:iss} ensure that the trajectory remains in $B_H(0,\varepsilon)$ throughout $[0,t]$ on $E_t$, and hence $L_t\in F_\varepsilon$.
Lemma~\ref{lem:wiener}$(ii)$ yields $\Prob_0(L_t\in F_\varepsilon)\ge c_1e^{-c_2t}$ for $t\ge1$.
Since $\delta_0\in\Lam(\gamma,M)$, the large deviation upper bound gives
\[
-\inf_{F_\varepsilon}I\ge\limsup_{t\to\infty}\frac1t\log\Prob_0(L_t\in F_\varepsilon)\ge-c_2 .
\]
As $I$ is a good rate function, the non-empty set $\{I\le c_2+1\}\cap F_\varepsilon$ is compact, and $I$ attains its minimum on this set at some $\sigma_\varepsilon$.
Although the chosen $M$ depends on $\mu_*$, Theorem~\ref{thm:A}$(iv)$ makes $I$ independent of $\gamma,M$.
The conclusions therefore concern the same rate function for every initial class.
\end{proof}
Finally, uniqueness of the equilibrium state for the zero potential identifies the zero of $I$.
\begin{proposition}\label{prop:zero}
Assume (A1)-(A3), and let $\sigma\in\cP(H)$. Then $I(\sigma)=0$ if and only if $\sigma=\mu_*$.
\end{proposition}

\begin{proof}
Proposition~\ref{prop:nontriv}$(i)$ gives $I(\mu_*)=0$.
Any other zero is an equilibrium state for $V=0$, so Lemma~\ref{lem:equil} identifies it with $h_0\mu_0$.
The proof of Proposition~\ref{prop:nontriv}$(i)$ gives $\lambda_0=1$ and $\mu_0=\mu_*$.
Since $P_t\mathbf1=\mathbf1$ and $\ip{\mathbf1}{\mu_*}=1$, the eigenfunction uniqueness in Theorem~\ref{thm:B}$(i)$ gives $h_0=1$.
Thus $\sigma=\mu_*$, proving uniqueness of the zero.
\end{proof}

\begin{proof}[Proof of Theorem~\ref{thm:A}]
Assertions $(i)$ and $(iv)$ follow from the preceding application of the Kifer-type criterion.
Proposition~\ref{prop:zero} proves $(ii)$.
For each $\varepsilon\in(0,1]$, Proposition~\ref{prop:nontriv}$(ii)$ provides $\sigma_\varepsilon$ supported in $B_H(0,\varepsilon)$ with finite rate and $\sigma_\varepsilon\ne\mu_*$.
Assertion $(ii)$ then gives $I(\sigma_\varepsilon)>0$, proving $(iii)$.
\end{proof}

\begin{proof}[Proof of Corollary~\ref{cor:level1}]
For $f\in C_b(H)$, the map $\Pi_f(\sigma):=\ip f\sigma$ is weakly continuous and $\Pi_f(L_t)=t^{-1}\int_0^tf(U_s)\dd s$.
The assertion therefore follows from Theorem~\ref{thm:A}$(i)$ by the contraction principle~\cite[Theorem~4.2.1]{DZ98}, with uniformity as in the remark after Corollary~1.4 of~\cite{JNPS18}.
\end{proof}

\subsection*{Declaration of competing interest}
The authors declare that they have no known competing financial interests or personal relationships that could have appeared to influence the work reported in this paper.

\subsection*{Data availability}
No datasets were generated or analysed during this study.

\subsection*{Funding}
Jianhua Huang is supported in part
by National Natural Science Foundation of China (No.12371198).
Jiangwei Zhang is supported by the China Postdoctoral Science Foundation (No.2026M794838).

\subsection*{Author Contributions}

All authors contributed equally to the conception, analysis, and preparation of this manuscript. All authors have read and approved the final manuscript.

\setlength{\bibsep}{4pt}


\begin{thebibliography}{99}

\bibitem{BGN23} P.-M. Boulvard, P. Gao, V. Nersesyan, Controllability and ergodicity of three dimensional primitive equations driven by a finite-dimensional force, Arch. Ration. Mech. Anal. 247 (2023) 2. \url{https://doi.org/10.1007/s00205-022-01835-8}

\bibitem{Cha06} D. Chae, Global regularity for the 2D Boussinesq equations with partial viscosity terms, Adv. Math. 203 (2006) 497-513. \url{https://doi.org/10.1016/j.aim.2005.05.001}

\bibitem{DM09} J. Duan, A. Millet, Large deviations for the Boussinesq equations under random influences, Stochastic Process. Appl. 119 (2009) 2052-2081. \url{https://doi.org/10.1016/j.spa.2008.10.004}

\bibitem{DZ14} G. Da Prato, J. Zabczyk, Stochastic Equations in Infinite Dimensions, second ed., Cambridge University Press, Cambridge, 2014.

\bibitem{DZ98} A. Dembo, O. Zeitouni, Large Deviations Techniques and Applications, second ed., Springer-Verlag, New York, 1998.

\bibitem{DP18} Z. Dong, X. Peng, Ergodicity of the 2D Navier-Stokes equations with degenerate multiplicative noise, Acta Math. Appl. Sin. Engl. Ser. 34 (2018) 97-118. \url{https://doi.org/10.1007/s10255-018-0728-z}

\bibitem{DV75} M.D. Donsker, S.R.S. Varadhan, Asymptotic evaluation of certain Markov process expectations for large time. I, Comm. Pure Appl. Math. 28 (1975) 1-47, \url{https://doi.org/10.1002/cpa.3160280102}; II, Comm. Pure Appl. Math. 28 (1975) 279-301, \url{https://doi.org/10.1002/cpa.3160280206}

\bibitem{FK06} J. Feng, T.G. Kurtz, Large Deviations for Stochastic Processes, American Mathematical Society, Providence, RI, 2006.

\bibitem{FGRT15} J. F\"oldes, N. Glatt-Holtz, G. Richards, E. Thomann, Ergodic and mixing properties of the Boussinesq equations with a degenerate random forcing, J. Funct. Anal. 269 (2015) 2427-2504. \url{https://doi.org/10.1016/j.jfa.2015.05.014}

\bibitem{FGRW16} J. F\"oldes, N.E. Glatt-Holtz, G. Richards, J.P. Whitehead, Ergodicity in randomly forced Rayleigh-B\'enard convection, Nonlinearity 29 (2016) 3309-3345. \url{https://doi.org/10.1088/0951-7715/29/11/3309}

\bibitem{Gou07a} M. Gourcy, A large deviation principle for 2D stochastic Navier-Stokes equation, Stochastic Process. Appl. 117 (2007) 904-927. \url{https://doi.org/10.1016/j.spa.2006.11.001}

\bibitem{Gou07b} M. Gourcy, Large deviation principle of occupation measure for stochastic Burgers equation, Ann. Inst. Henri Poincar\'e Probab. Stat. 43 (2007) 441-459. \url{https://doi.org/10.1016/j.anihpb.2006.07.003}

\bibitem{GL21} V.M. G\'alfi, V. Lucarini, Fingerprinting heatwaves and cold spells and assessing their response to climate change using large deviation theory, Phys. Rev. Lett. 127 (2021) 058701. \url{https://doi.org/10.1103/PhysRevLett.127.058701}

\bibitem{GHM18} N.E. Glatt-Holtz, D.P. Herzog, J.C. Mattingly, Scaling and saturation in infinite-dimensional control problems with applications to stochastic partial differential equations, Ann. PDE 4 (2018), Paper No. 16. \url{https://doi.org/10.1007/s40818-018-0052-1}

\bibitem{HL05} T.Y. Hou, C. Li, Global well-posedness of the viscous Boussinesq equations, Discrete Contin. Dyn. Syst. 12 (2005) 1-12.

\bibitem{HM06} M. Hairer, J.C. Mattingly, Ergodicity of the 2D Navier-Stokes equations with degenerate stochastic forcing, Ann. of Math. (2) 164 (2006) 993-1032. \url{https://doi.org/10.4007/annals.2006.164.993}

\bibitem{JNPS15} V. Jak\v si\'c, V. Nersesyan, C.-A. Pillet, A. Shirikyan, Large deviations from a stationary measure for a class of dissipative PDEs with random kicks, Comm. Pure Appl. Math. 68 (2015) 2108-2143. \url{https://doi.org/10.1002/cpa.21568}

\bibitem{JNPS18} V. Jak\v si\'c, V. Nersesyan, C.-A. Pillet, A. Shirikyan, Large deviations and mixing for dissipative PDEs with unbounded random kicks, Nonlinearity 31 (2018) 540-596. \url{https://doi.org/10.1088/1361-6544/aa99a7}

\bibitem{Kif90} Y. Kifer, Large deviations in dynamical systems and stochastic processes, Trans. Amer. Math. Soc. 321 (1990) 505-524. \url{https://doi.org/10.1090/S0002-9947-1990-1025756-7}

\bibitem{KM23a} A. Kumar, M.T. Mohan, Large deviation principle for occupation measures of two dimensional stochastic convective Brinkman-Forchheimer equations, Stoch. Anal. Appl. 41 (2023) 214-256. \url{https://doi.org/10.1080/07362994.2021.2005626}

\bibitem{KM23b} A. Kumar, M.T. Mohan, Large deviation principle for occupation measures of stochastic generalized Burgers-Huxley equation, J. Theoret. Probab. 36 (2023) 661-709. \url{https://doi.org/10.1007/s10959-022-01180-2}

\bibitem{KS12} S. Kuksin, A. Shirikyan, Mathematics of Two-Dimensional Turbulence, Cambridge University Press, Cambridge, 2012.

\bibitem{Maj03} A. Majda, Introduction to PDEs and Waves for the Atmosphere and Ocean, Courant Lecture Notes in Mathematics, vol. 9, American Mathematical Society, Providence, RI, 2003.

\bibitem{MN18a} D. Martirosyan, V. Nersesyan, Local large deviations principle for occupation measures of the stochastic damped nonlinear wave equation, Ann. Inst. Henri Poincar\'e Probab. Stat. 54 (2018) 2002-2041. \url{https://doi.org/10.1214/17-AIHP863}

\bibitem{Ner19} V. Nersesyan, Large deviations for the Navier-Stokes equations driven by a white-in-time noise, Ann. H. Lebesgue 2 (2019) 481-513. \url{https://doi.org/10.5802/ahl.23}

\bibitem{NPX23} V. Nersesyan, X. Peng, L. Xu, Large deviations principle via Malliavin calculus for the Navier-Stokes system driven by a degenerate white-in-time noise, J. Differential Equations 362 (2023) 230-249. \url{https://doi.org/10.1016/j.jde.2023.03.004}

\bibitem{LL26} R. Liu, K. Lu, Large deviations for 2D Navier-Stokes equations driven by a periodic force and a degenerate noise, Ann. Probab. 54 (2026) 705-767. \url{https://doi.org/10.1214/25-AOP1777}

\bibitem{PHZ20a} X. Peng, J. Huang, Y. Zheng, Exponential mixing for the fractional magneto-hydrodynamic equations with degenerate stochastic forcing, Commun. Pure Appl. Anal. 19 (2020) 4479-4506. \url{https://doi.org/10.3934/cpaa.2020204}

\bibitem{PHZ20b} X. Peng, J. Huang, R. Zhang, Ergodicity and exponential mixing of the real Ginzburg-Landau equation with a degenerate noise, J. Differential Equations 269 (2020) 3686-3720. \url{https://doi.org/10.1016/j.jde.2020.03.013}

\bibitem{PX26} X. Peng, L. Xu, Large deviations principle for 2D Navier-Stokes equations with space-time localized noise, Commun. Math. Stat. (2026) 1-17. \url{https://doi.org/10.1007/s40304-025-00478-x}

\bibitem{NR25} V. Nersesyan, M. Rissel, Global controllability of Boussinesq flows by using only a temperature control, Arch. Ration. Mech. Anal. 249 (2025), Article 58. \url{https://doi.org/10.1007/s00205-025-02128-6}

\bibitem{Nua06} D. Nualart, The Malliavin Calculus and Related Topics, second ed., Springer-Verlag, Berlin, 2006.

\bibitem{Ris25} M. Rissel, Controllability of Boussinesq flows driven by finite-dimensional and physically localized forces, arXiv:2506.19764 (2025).

\bibitem{Ris26} M. Rissel, Approximate controllability of Boussinesq channel flows only through the temperature, SIAM J. Control Optim. 64 (2026) 59-75. \url{https://doi.org/10.1137/25M1746276}

\bibitem{RWB18} F. Ragone, J. Wouters, F. Bouchet, Computation of extreme heat waves in climate models using a large deviation algorithm, Proc. Natl. Acad. Sci. USA 115 (2018) 24-29. \url{https://doi.org/10.1073/pnas.1712645115}

\bibitem{Tyc35} A. Tychonoff, Ein Fixpunktsatz, Math. Ann. 111 (1935) 767-776. \url{https://doi.org/10.1007/BF01472256}

\bibitem{WXX21} R. Wang, J. Xiong, L. Xu, Large deviation principle of occupation measures for non-linear monotone SPDEs, Sci. China Math. 64 (2021) 799-822. \url{https://doi.org/10.1007/s11425-018-9482-1}

\bibitem{Wu01} L. Wu, Large and moderate deviations and exponential convergence for stochastic damping Hamiltonian systems, Stochastic Process. Appl. 91 (2001) 205-238. \url{https://doi.org/10.1016/S0304-4149(00)00061-2}

\bibitem{Zhao25} M. Zhao, Level-3 large deviations for the white-forced 2D Navier-Stokes system in a bounded domain, arXiv:2506.14119 (2025).

\bibitem{ZH16} Y. Zheng, J. Huang, Large deviation principle for stochastic Boussinesq equations driven by L\'evy noise, J. Math. Anal. Appl. 439 (2016) 523-550. \url{https://doi.org/10.1016/j.jmaa.2016.02.054}

\end{thebibliography}
\end{document}